\documentclass[11pt]{article}
\usepackage[english]{babel}

\usepackage{amsmath,amssymb,amsthm}
\usepackage{bbm}
\usepackage{mathtools}

\usepackage{pgfplots}
\usepackage{tikz}
\usepackage{multirow}
\usepackage{subcaption}
\usepackage[scale=0.7]{geometry}

\pgfplotsset{
	log y ticks with fixed point/.style={
		yticklabel={
			\pgfkeys{/pgf/fpu=true}
			\pgfmathparse{exp(\tick)}%
			\pgfmathprintnumber[fixed relative, precision=3]{\pgfmathresult}
			\pgfkeys{/pgf/fpu=false}
		},
		y tick label style={/pgf/number format/1000 sep=},
	}
}

\makeatletter

\newcommand{\printslope}[4]{
	\tikzset{fixed point arithmetic}
	\def\nero@printslope@orderlist{#1}
	\edef\nero@printslope@xpos{#2}
	\edef\nero@printslope@ypos{#3}
	\edef\nero@printslope@width{#4}
	\pgfmathparse{\nero@printslope@xpos+\nero@printslope@width}
	\edef\nero@printslope@px{\pgfmathresult}
	\edef\nero@printslope@py{\nero@printslope@ypos}
	\edef\nero@printslope@qx{\nero@printslope@xpos}
	\edef\nero@printslope@ry{\nero@printslope@ypos}
	\foreach \nero@printslope@order in {#1}{
		\pgfmathparse{
			((\nero@printslope@px/\nero@printslope@xpos)^(\nero@printslope@order))*\nero@printslope@ypos}
		\edef\nero@printslope@qy{\pgfmathresult}
		\edef\nero@aux1{\noexpand\draw[line width=0.6pt]
			(axis cs:\nero@printslope@xpos,\nero@printslope@ry)
			-- (axis cs:\nero@printslope@qx,\nero@printslope@qy)
			-- (axis cs:\nero@printslope@px,\nero@printslope@py);}
		\nero@aux1
		\pgfmathparse{10^((ln(\nero@printslope@ry)+ln(\nero@printslope@qy))/(ln(10)*2))}
		\edef\nero@printslope@labelpos{\pgfmathresult}
		\edef\nero@aux2{\noexpand\node[anchor=east] at
			(axis cs:\nero@printslope@qx,\nero@printslope@labelpos)
			{\noexpand\small \nero@printslope@order};}
		\nero@aux2
		\global\edef\nero@printslope@ry{\nero@printslope@qy}
	}
	\draw[line width=0.6pt] (axis cs:\nero@printslope@xpos,\nero@printslope@ypos)
	|- (axis cs:\nero@printslope@px,\nero@printslope@py);
	\pgfmathparse{10^((ln(\nero@printslope@px)+ln(\nero@printslope@xpos))/(ln(10)*2))}
	\edef\nero@printslope@labelpos{\pgfmathresult}
	\node[anchor=north] at (axis cs:\nero@printslope@labelpos,\nero@printslope@ypos) {\small 1};
}

\makeatother

\newtheorem{thm}{Theorem}
\newtheorem{lem}[thm]{Lemma}
\newtheorem{pf}{Proof}
	
\title{
	\bf{An all speed second order well-balanced IMEX relaxation scheme for the Euler equations with gravity}
}
\author{
	Andrea Thomann\footnote{Dipartimento di Scienze e Alta Tecnologia, Universit\`a degli Studi dell'Insubria, Via Valleggio 11, 22100 Como, Italy}  \footnote{Correspondence to: Andrea Thomann, University of Insubria, Via Valleggio 11, 22100 Como (CO), Italy, Email: acthomann@uninsubria.it} \footnote{Marie Sk\l odowska-Curie fellow of the Istituto Nazionale di Alta Matematica Francesco Severi, Rome, Italy},
Gabriella Puppo\footnote{Dipartimento di Matematica, La Sapienza Universit\`a di Roma, Piazzale Aldo Moro 5, 00185 Roma, Italy},
Christian Klingenberg\footnote{Fakult\"at f\"ur Mathematik und Informatik, Universit\"at W\"urzburg, Emil-Fischer-Str. 40, 97074 W\"urzburg, Germany}
}
\date{\today}

\begin{document}
\maketitle

\section*{Abstract}
We present an implicit-explicit well-balanced finite volume scheme for the Euler equations with a gravitational source term which is able to deal also with low Mach flows. 
To visualize the different scales we use the non-dimensionalized equations on which we apply a pressure splitting and a Suliciu relaxation. 
On the resulting model, we apply a splitting of the flux into a linear implicit and an non-linear explicit part that leads to a scale independent time-step. 
The explicit step consists of a Godunov type method based on an approximative Riemann solver where the source term is included in the flux formulation.
We develop the method for a first order scheme and give an extension to second order. 
Both schemes are designed to be well-balanced, preserve the positivity of density and internal energy and have a scale independent diffusion. 
We give the low Mach limit equations for well-prepared data and show that the scheme is asymptotic preserving.
These properties are numerically validated by various test cases. 

\paragraph{Keywords} IMEX scheme, Suliciu relaxation, Euler equations with gravity, non-dimensional, well-balanced, positivity preserving, asymptotic preserving

\section{Introduction}
The aim of this paper is the construction of an all speed scheme for the Euler equations of gas dynamics with a given gravitational source term in multiple space dimensions. 
Applications of this model can be found for example in astrophysics and meteorology. 
A broad overview is given in the review of Klein \cite{Klein2010} where it is demonstrated that
atmospheric flows can have large scale differences.
To reflect those scales in the equations, we use the non-dimensionalised version which is characterized by the reference Mach and Froude numbers denoted by $M$ and $Fr$ respectively.

In the homogeneous case the behaviour of the fluid changes depending on the Mach number only. 
It ranges from compressible flow for large Mach numbers to the incompressible limit equations for $M$ going to zero. 
The derivation of the limit equations can be found eg. in  \cite{KlainermanMajda1982,Dellacherie2010,Schochet2005} and references therein.
To accurately approximate all speed flows, asymptotic preserving (AP) schemes are well suited since they are consistent with the limit behaviour as $M$ tends to zero.
The development of those schemes is an active field of research and we refer to the review of Jin \cite{Jin2012} for an introduction and \cite{BarsukowHohmKBRoe2019} for a recent work on an active flux method for linear acoustics.
An important role in the achievement of the AP property is played by the splitting of the pressure following the studies of Klein \cite{Klein1995,KleinEtAl2001} as used in the schemes \cite{NoelleEtAl2014weakly,CordierDegondKumbaro2012,ThomannZenkPuppoKB2019,BerthonKlingenbergZenk2018}.
In \cite{ThomannZenkPuppoKB2019,BerthonKlingenbergZenk2018} the pressure splitting is combined with a Suliciu relaxation approach \cite{Suliciu1990} which allows for an easy construction of Riemann Solvers.
An example for a Jin-Xin relaxation approach \cite{JinXin1995} can be found in \cite{AbbateIolloPuppo2017}.

Since for explicit schemes the time step is restricted by the inverse of the largest wave speed which scales with $1/M$, explicit schemes are not practical for low Mach applications.
Therefore implicit \cite{BerthonKlingenbergZenk2018,AbbateIolloPuppo2017} or implicit-explicit (IMEX) schemes \cite{CordierDegondKumbaro2012,BoscarinoRussoScandurra2018,DimarcoEtAl2018} are used to have Mach number independent time step.

The presence of the source term makes it interesting to look at steady states.
For zero velocity, we find the hydrostatic equilibrium, that is characterized as the balance of the pressure gradient with the weight of the fluid. 
Most atmospheric-flow phenomena may be understood as perturbations of such a balanced background state.
The scope of well-balanced schemes is to maintain the background atmosphere at machine precision to be able to resolve those small perturbations accurately.
Since the shape of the equilibrium state depends on the underlying pressure law there are schemes focused on well-balancing a specific class of equilibria, for example isothermal and polytropic atmospheres \cite{DesveauxZenkBerthonKlingenberg2016} or equilibria with constant entropy \cite{KaeppeliMishra2014}.
The latter was extended in \cite{KaeppeliMishra2016} to the preservation of hydrostatic equilibria with arbitrary entropy stratification using a second order reconstruction of the discrete equilibrium equation. 
A different approach can be found in \cite{GaburroCastroDumbser2018}, where the well-balanced property is achieved by using path-conservative finite volumes schemes. 
Higher order well-balanced schemes can be realized by using a high order hydrostatic reconstruction, as done in \cite{XingShu2013,KlingenbergPuppoSemplice2019,GrosheintzKaeppeli2019}. 
Since our aim is to exactly well-balance arbitrary hydrostatic equilibria, we follow the approach used in
\cite{KlingenbergPuppoSemplice2019,GhoshConstantinescu2015,ThomannZenkKlingenberg2019} and rewrite the gravitational potential in terms of a reference equilibrium state.
Note that the above mentioned well-balancing techniques were developed for the compressible regime. 
To have a well-balanced scheme that is applicable in the low Mach, low Froude regime, we extend the second order AP IMEX scheme developed for the homogeneous Euler equations \cite{ThomannZenkPuppoKB2019} 
to include also a gravitational source term. The new scheme is designed to inherit the nice properties of the homogeneous case. In particular, it preserves the positivity of the density and of the pressure, enjoys a Mach number independent numerical diffusion, and it can be easily extended to second order. 

To our knowledge, this is the first case in which the construction of a well balanced scheme for general equilibria is addressed which, at the same time, preserves asymptotic properties in the low Mach regime under a gravitational field for the full Euler equations.
%  \color{red} Preliminary results markus thesis on well-balanced in the low Mach for isothermal and polytropic \color{black}
We show the AP property of the scheme by proving that it preserves the divergence free constraint in the zero Mach number limit when starting from well prepared initial data. 
The limit equations are given by the incompressible Euler equations in a gravitational field. 
% to a gravitational source term and inherit the nice properties as preserving the positivity and a natural extension to second order.
% Therefore we have to modify, following \cite{BerthonKlingenbergZenk2018,ThomannZenkKlingenberg2019}, the Suliciu relaxation approach to incorporate the source term and at the same time keep the ordered wave structure which guarantees an easy construction of a Riemann solver. 
% As in the homogeneous case, we split the pressure into a slow and a fast component, as discussed in \cite{Klein1995}.  
%To guarantee a scale independent time step, the flux function is split into a linear implicit part only involving relaxation variables and an explicit part which contains the original non-linear flux function of the Euler equations and concerns the update of the physical variables.   
%Both splittings are essential to ensure the AP property and the scale independent diffusion.
%To show the AP property, we define a set of well-prepared data which consists at leading order of the hydrostatic equilibrium.
% as well as the divergence free property of the velocity and its orthogonality to the direction of the gravitation. 
%In the limit this results in the incompressible Euler equations with a gravitational source term. 
Similar results were found in \cite{BispenLukacovaYelash2017} for the isentropic case with potential temperature. 
We refer to \cite{FeireislKBKremlMarkfelder2017Arxiv,FeireislKlingenbergMarkfelder2019} for theoretical studies on the isothermal and isentropic case with a one component linear gravitational field and to \cite{BarsukovEdelmannKBRoepke2017} for a low Mach scheme that allows for a gravitational source term, but lacks the well-balanced property.

The paper is organized as follows. 
In Section 2, we introduce the equations, notion of hydrostatic equilibria and the limit equations. 
Then we give the derivation of the Suliciu type relaxation model in Section 3. 
The time semi-discrete scheme with the flux splitting together with the Mach number expansion of the fast pressure and the asymptotic preserving property are discussed in Section 4.
Subsequent, we give the derivation of the fully discrete scheme which includes a Godunov type finite volume scheme based on an approximative Riemann solver in the explicit part. 
We show that the scheme is well-balanced and that it preserves the positivity of the density and internal energy.
The section ends with the extension to second order. 
All properties are numerically validated in Section 6. 
In particular, we give an example of low Mach flow, starting from well prepared initial data, and we study a low Mach stationary vortex in a gravitational field, with a test we derived from the classical Gresho vortex test case \cite{MiczekRoepkeEdelmann2015}.
We conclude the numerical tests with a simulation of a rising hot air bubble which arises in meteorology. 

A section of conclusion completes this paper. 

\section{The Euler equations with a gravitational source term}
% gamma is the adiabatic constant
The Euler equations with a gravitational source term in $d$ dimensions are given by 
\begin{align}
\begin{split}
\rho_t + \nabla \cdot (\rho \mathbf{u}) &= 0,\\
(\rho \mathbf{u})_t + \nabla \cdot \left(\rho \mathbf{u}\otimes\mathbf{u}\right) + \nabla p &= -\rho \nabla \Phi,\\
E_t + \nabla \cdot ( \mathbf{u} (E + p)) &= -\rho \mathbf{u} \cdot \nabla \Phi\\
\end{split}
\label{sys:Euler_G}
\end{align}
where the total energy $E$ is given by 
\begin{equation*}
\label{eq:Energy_dim}
E = \rho e + \frac{1}{2} \rho |\mathbf{u}|^2.
\end{equation*}
Here, $\rho$ denotes the density, $\mathbf{u} \in \mathbb{R}^d$ the velocity vector, $e$ the internal energy and $\Phi: \mathbb{R}^d \to \mathbb{R}$ is a given smooth stationary potential.
The pressure is given by the ideal gas law $ p = (\gamma -1) \rho e$ and the speed of sound is denoted by $c$.

To make the impact of slow and fast scales evident in the equations, we rewrite \eqref{sys:Euler_G} in its non-dimensional form by decomposing all variables $\varphi$ into a scalar reference value $\varphi_r$, that contains the units, and a non-dimensional quantity $\tilde \varphi$: 
\begin{equation}
\label{eq:decomposition}
\varphi = \varphi_r ~ \tilde{\varphi}.
\end{equation}
Choosing the reference length $x_r$, time $t_r$, density $\rho_r$, sound speed $c_r$ and gravitational acceleration $\Phi_r$, we can compute the missing reference values as
\begin{equation}
\label{eq:relation_refValues}
u_r = \frac{x_r}{t_r}, ~~\text{ and }~~p_r = \rho_r c_r^2.
\end{equation}
Inserting the decomposition \eqref{eq:decomposition} in the dimensional equations \eqref{sys:Euler_G} and using the relations \eqref{eq:relation_refValues}, we arrive at the non-dimensional Euler equations with a gravitational source term:
\begin{align}
\begin{split}
\rho_t + \nabla \cdot (\rho \mathbf{u}) &= 0\\
(\rho \mathbf{u})_t + \nabla \cdot \left(\rho \mathbf{u}\otimes\mathbf{u}\right) + \frac{1}{M^2} \nabla p &= -\frac{1}{Fr^2}\rho \nabla \Phi\\
E_t + \nabla \cdot ( \mathbf{u} (E + p)) &= -\frac{M^2}{Fr^2}\rho \mathbf{u} \cdot \nabla \Phi.
\end{split}
\label{sys:Euler_LMG}
\end{align}
For simplicity, we have dropped the tilde and, if not otherwise mentioned, we will use the non-dimensional variables throughout this paper.
The total energy of system \eqref{sys:Euler_LMG} is given by
\begin{equation*}
%\label{eq:Etot}
E = \rho e + \frac{1}{2}M^2 \rho |\mathbf{u}|^2.
\end{equation*}
Equations \eqref{sys:Euler_LMG} depend on two non-dimensional quantities, the Mach number $M$ and the Froude number $Fr$. 
The Mach number is defined as the ratio between the velocity of the gas and the sound speed
\begin{equation*}
M = \frac{u_r}{c_r}
\end{equation*}
and the Froude number is defined as the ratio between the velocity of the gas and the velocity introduced by the gravitational acceleration
\begin{equation*}
Fr = \frac{u_r}{\sqrt{\Phi_r}}.
\end{equation*}

\subsection{Hydrostatic equilibria}

Hydrostatic equilibria are stationary solutions of \eqref{sys:Euler_LMG} that satisfy
\begin{align}
\label{eq:hydro_steady}
\begin{split}
\mathbf{u} &= 0, \\
\frac{1}{M^2}\nabla p &= - \frac{1}{Fr^2}\rho \nabla \Phi.
\end{split}
\end{align}
Solutions to \eqref{eq:hydro_steady} are not unique and depending on the relation between the pressure and the density they can have completely different behaviour. 
To demonstrate this, let us for a moment consider the following class 
\begin{equation}
\label{eq:EOS}
p = \chi \rho^\Gamma
\end{equation}
with constants $\chi > 0$, $\Gamma \in (0,\infty)$. 
For the class of equation of states \eqref{eq:EOS}, we obtain for $\Gamma = 1$ (isothermal) with a constant $c \in \mathbb{R}$ and $\chi = RT$
\begin{align}
\label{eq:isoth}
\begin{cases}
\mathbf{u}(\mathbf{x}) &= 0, \\
\overline \rho(\mathbf{x}) &= \exp\left(\frac{c - \frac{M^2}{Fr^2}\Phi(\mathbf{x})}{RT}\right)\\
\overline p(\mathbf{x}) &= RT \exp\left(\frac{c - \frac{M^2}{Fr^2}\Phi(\mathbf{x})}{RT}\right)
\end{cases}
\end{align}
and for $\Gamma \in (0,1) \cup (1,\infty)$ (polytropic) with a constant $c \in \mathbb{R}$
\begin{align}
\label{eq:poly}
\begin{cases}
\mathbf{u}(\mathbf{x}) &= 0, \\
\overline \rho(\mathbf{x}) &= \left(\frac{\Gamma -1}{\chi \Gamma} (c - \frac{M^2}{Fr^2}\Phi(\mathbf{x}))\right)^{\frac{1}{\Gamma -1}}\\
\overline p(\mathbf{x}) &= \chi^{\frac{1}{1 - \Gamma}}\left(\frac{\Gamma -1}{\chi \Gamma} (c - \frac{M^2}{Fr^2}\Phi(\mathbf{x}))\right)^{\frac{\Gamma}{\Gamma -1}}.
\end{cases}
\end{align}

Since arbitrary solutions $\overline \rho$ and $\overline p$ of the hydrostatic equilibrium \eqref{eq:hydro_steady} are stationary, we follow \cite{GhoshConstantinescu2015} and define two time-independent positive functions $\alpha(\mathbf{x}) = \overline \rho(\mathbf{x})$ and $\beta(\mathbf{x}) = \overline p(\mathbf{x})$ representing the equilibrium density and pressure respectively. 
Since $\alpha,\beta$ satisfy \eqref{eq:hydro_steady}, we can find a new relation for $\nabla \Phi$ due to the following equivalent description
\begin{equation}
\label{eq:DefBeta}
\frac{1}{M^2}\nabla \beta = - \frac{1}{Fr^2}\alpha \nabla \Phi ~~~\Leftrightarrow~~~\nabla \Phi = - \frac{Fr^2}{M^2} \frac{\nabla \beta}{\alpha}.
\end{equation}
With this definition of the gravitational potential, we can rewrite \eqref{sys:Euler_LMG} into
\begin{align}
\begin{split}
\rho_t + \nabla \cdot (\rho u) &= 0,\\
(\rho \mathbf{u})_t + \nabla \cdot \left(\rho \mathbf{u}\otimes\mathbf{u}\right) + \frac{1}{M^2} \nabla p &= \frac{1}{M^2}\frac{\rho}{\alpha} \nabla \beta,\\
E_t + \nabla \cdot ( \mathbf{u} (E + p)) &= \frac{\rho}{\alpha}\mathbf{u} \cdot \nabla \beta.
\end{split}
\label{sys:Euler_LMG_mod}
\end{align}
We emphasize, that the reference equilibrium has to be known in advance. 
In general, this is not a restriction, because in many applications the equilibrium solutuion of interest is known in advance.
Note, that the equations \eqref{sys:Euler_LMG_mod} are only depending on the Mach number, but the dependence on the Froude number is implicitly given in the definition of $\beta$ in \eqref{eq:DefBeta}.

\subsection{The low Mach limit} 
%To have the balance of pressure gradient and source term at leading order results in a scaling of $\mathcal{O}(M^2) = \mathcal{O}(Fr^2)$.
To model perturbations of an equilibrium state, we assume in accordance with
\cite{BispenLukacovaYelash2017,FeireislKBKremlMarkfelder2017Arxiv,FeireislKlingenbergMarkfelder2019} that $\mathcal{O}(Fr^2) = \mathcal{O}(M^2)$.
To analyse multi-scale effects and the formal asymptotic behaviour of \eqref{sys:Euler_LMG}, we express the variables in form of a Mach number expansion and compare the orders of terms in $M$.
The expansions are given by
\begin{align}
\begin{array}{cc}
\rho = \rho_0 + M \rho_1 + M^2 \rho_2 + \mathcal{O}(M^3), & \mathbf{u} = \mathbf{u}_0 + M \mathbf{u}_1 + M^2 \mathbf{u}_2 + \mathcal{O}(M^3),\\ 
e = e_0 + M e_1 + M^2 e_2 + \mathcal{O}(M^3), & p = p_0 + M p_1 + M^2 p_2 + \mathcal{O}(M^3).
\end{array}
\label{exp:MFr}
\end{align}

Inserting the expansion \eqref{exp:MFr} into the Euler equations \eqref{sys:Euler_LMG_mod} and collecting the terms of order $\mathcal{O}(M^{-2})$, we have
\begin{align}
\nabla p_{0} = - \rho_0 \nabla \Phi.
\label{eq:WPM-2}
\end{align}
For the $\mathcal{O}(M^{-1})$ terms, we find
\begin{align}
\nabla p_{1} = - \rho_1 \nabla \Phi.
\label{eq:WPM-1}
\end{align}
This means that the couples $p_0, \rho_0$ and $p_1, \rho_1$ fulfil the hydrostatic equilibrium and thus are time-independent. 
Using this in the $\mathcal{O}(M^0)$ terms, we obtain
\begin{align*}
\begin{split}
\nabla \cdot (\rho_0 \mathbf{u}_0) &= 0, \\
\partial_t \mathbf{u}_0 + \mathbf{u}_0 \cdot \nabla \mathbf{u}_0 + \frac{\nabla p_2}{\rho_0} &= - \frac{\rho_2 \nabla \Phi}{\rho_0}, \\
\nabla \cdot \mathbf{u}_0 &= \frac{\mathbf{u}_0 \cdot \nabla \Phi}{c_0^2},
\end{split}
\end{align*}
where we have used $c_0^2 = \gamma \frac{p_0}{\rho_0}$.
We define the set of well-prepared data for a given potential $\Phi$ as 
\begin{align}
\begin{split}
\Omega_{wp} = \left\lbrace w \in \mathbb{R}^{d+2}\left|\nabla p_0 = -\rho_0 \nabla \Phi, ~\nabla p_1 \right.\right.&= -\rho_1 \nabla \Phi, ~\nabla \cdot (\rho_0 \mathbf{u}_0) = 0,\\
&\hskip-3cm\left.\nabla \cdot \mathbf{u}_0 = 0, ~\mathbf{u}_0 \cdot \nabla \Phi = 0 \right\rbrace.
\end{split}
\label{eq:Omegawp}
\end{align}
Analogously, we define the well-prepared data for given $\alpha,\beta$ for the modified equations \eqref{sys:Euler_LMG_mod}
\begin{align}
\begin{split}
\Omega_{wp}^{\alpha\beta} = \left\lbrace w \in \mathbb{R}^{d+2}\left| \nabla p_0 = \rho_0 \frac{\nabla \beta}{\alpha}, ~\nabla p_1 \right.\right.&= \rho_1 \frac{\nabla \beta}{\alpha},~\nabla \cdot (\rho_0 \mathbf{u}_0) = 0, \\
&\hskip-3cm\left.\nabla \cdot \mathbf{u}_0 = 0, ~ \mathbf{u}_0 \cdot \frac{\nabla \beta}{\alpha} = 0 \right\rbrace.
\end{split}
\label{eq:Omegawpab}
\end{align}
This means that the pressure and density fulfil the hydrostatic equilibrium up to a perturbation of $M^2$, the the first component of the velocity field is divergence free and orthogonal to $\nabla \Phi$.
Thus we obtain as the limit equations, the incompressible Euler equations with a gravitational source term 
\begin{align}
\label{sys:LimitEQ}
\begin{split}
\nabla \cdot (\rho_0 \mathbf{u}_0) &= 0,\\
\partial_t \mathbf{u}_0 + \mathbf{u}_0 \cdot \nabla \mathbf{u}_0 + \frac{\nabla p_2}{\rho_0} &= - \frac{\rho_2 \nabla \Phi}{\rho_0}, \\
\nabla \cdot \mathbf{u}_0 = 0,  ~~\mathbf{u}_0 \cdot \nabla \Phi &= 0.
\end{split}
\end{align} 

\section{Suliciu Relaxation model}
Using a Suliciu Relaxation approach \cite{Suliciu1990,Bouchut2004,CoquelPerthame1998} is one way of simplifying the non-linear structure of the Euler equations \eqref{sys:Euler_G}. 
The derivation of the relaxation model follows the argument given in \cite{CordierDegondKumbaro2012,ThomannZenkPuppoKB2019,BerthonKlingenbergZenk2018}. 
In the spirit of Klein \cite{Klein1995}, we apply in the momentum and energy equation a splitting of the pressure $p$ into a slow and a fast component 
\begin{equation*}
\frac{p}{M^2} = p + \frac{1 - M^2}{M^2} p.
\end{equation*} 
The aim is to relax both the slow and the fast pressure in a Suliciu relaxation manner. 
The pressure in relaxation equilibrium is then characterized by 
\begin{equation*}
p = M^2 \pi + (1 - M^2) \psi,
\end{equation*} 
where $\pi$ is the approximation of the slow and $\psi$ of the fast part.
To obtain the evolution of $\pi$, we can directly apply the Suliciu relaxation technique which leads to the addition of the following equation in conservation form
\begin{equation*}
(\rho \pi)_t + \nabla \cdot \left(\rho \pi \mathbf{u}\right) + a^2 \nabla \cdot \mathbf{u} = \frac{\rho}{\varepsilon}( p - \pi).
\end{equation*} 
As discussed in \cite{BerthonKlingenbergZenk2018}, applying this Suliciu relaxation technique also on the fast pressure does not lead to scheme that is accurate for small Mach numbers. 
Instead a relaxation equation for the velocity $\mathbf{\hat{u}}$ coupled the pressure $\psi$ is added.  
We apply the same strategy as in the homogeneous case described in \cite{ThomannZenkPuppoKB2019,BerthonKlingenbergZenk2018}.
Here in additionh, the influence of the source term in the momentum equation has to be taken into account.
As a consequence, the source term will also appear in the relaxation equation for $\mathbf{\hat{u}}$.
The relaxation model is developed under the following objectives:
\begin{itemize}
	\item It has ordered eigenvalues that lead to a clear wave structure and make it especially easy to construct a Riemann solver.
	\item It is a stable diffusive approximation of the non-dimensional Euler equations with gravitational source term \eqref{sys:Euler_LMG_mod}.
	\item The resulting numerical scheme has Mach number independent diffusion.
\end{itemize}
The achievement of the first objective depends also on the treatment of the source term, since it is associated over $\beta$ with a 0 eigenvalue.
Following \cite{DesveauxZenkBerthonKlingenberg2016}, we remove the 0 eigenvalue by relaxing also $\beta$.
It is approximated by a new variable $Z$ that is transported with $\mathbf{u}$ as
\begin{equation*}
Z_t + \mathbf{u}\cdot\nabla Z = \frac{1}{\varepsilon}(\beta - Z).
\end{equation*}
This associates the source term with the eigenvalue $\mathbf{u}$.
Since the evolution of $\alpha$ is constant in time, we consider it as a given time independent function and will omit its evolution in the relaxation model.
All this considerations lead to the following relaxation model in conservation form:
\begin{align}
\begin{split}
\rho_t + \nabla \cdot (\rho u) &= 0\\
(\rho \mathbf{u})_t + \nabla \cdot \left(\rho \mathbf{u}\otimes\mathbf{u}\right) + \nabla \pi + \frac{1-M^2}{M^2} \nabla \psi &= \frac{1}{M^2}\frac{\rho}{\alpha} \nabla Z,\\
E_t + \nabla \cdot ( \mathbf{u} (E + M^2 \pi + (1 - M^2) \psi)) &= \frac{\rho}{\alpha} \mathbf{u} \cdot \nabla Z, \\
(\rho \pi)_t + \nabla \cdot (\rho \mathbf{u} \pi + a^2 \mathbf{u}) &= \frac{\rho}{\varepsilon}(p - \pi), \\
(\rho \mathbf{\hat{u}})_t + \nabla \cdot \left(\rho \mathbf{u} \otimes \mathbf{\hat{u}}\right) + \frac{1}{M^2} \nabla\psi &= \frac{1}{M^2}\frac{\rho}{\alpha} \nabla Z + \frac{\rho}{\varepsilon} (\mathbf{u} - \mathbf{\hat{u}}), \\
(\rho \psi)_t + \nabla \cdot (\rho \mathbf{u} \psi + a^2 \mathbf{\hat{u}}) &= \frac{\rho}{\varepsilon}(p - \psi),\\
(\rho Z)_t + \nabla (\rho \mathbf{u} Z) &= \frac{\rho}{\varepsilon}(\beta - Z).
\end{split}
\label{sys:full_relax_ELMG}
\end{align}
The following lemma sums up some properties of system \eqref{sys:full_relax_ELMG}.
\begin{lem}
	\label{lem:properties}
	The relaxation system \eqref{sys:full_relax_ELMG} is hyperbolic and is a stable diffusive approximation of \eqref{sys:Euler_LMG_mod} under the Mach number independent sub-characteristic condition for the relaxation parameter $a > \rho \sqrt{\partial_\rho p(\rho,e)}$.
	It has the following linearly degenerate eigenvalues 
	\begin{equation*}
	\lambda^u = \mathbf{u}, \lambda^\pm = \mathbf{u} \pm \frac{a}{\rho}, \lambda_M^\pm = \mathbf{u} \pm \frac{a}{M\rho}.
	\end{equation*} 
	%	For $M < 1$, the eigenvalues have the ordering
	%	\begin{equation*}
	%		\lambda_M^- < \lambda^- < \lambda^u < \lambda^+ < \lambda_M^+.
	%	\end{equation*}
\end{lem}
\begin{pf}
	The proof follows the same lines of  \cite{ThomannZenkPuppoKB2019,DesveauxZenkBerthonKlingenberg2016,ThomannZenkKlingenberg2019} and is omitted.
\end{pf}
\vskip0.25cm
Note, that in the case of $M=1$, the waves associated with $\lambda_M^\pm$ and $\lambda^\pm$ collapse to $\lambda^\pm$. 
For simplicity, we will refer to system \eqref{sys:Euler_LMG_mod} as
\begin{equation}
\label{sys:short_mod_EulerLMG}
w_t + \nabla \cdot f(w) = s(w).
\end{equation}
where $w = (\rho, \rho \mathbf{u}, E)^T$ denotes the vector of physical variables, while the flux function $f(w)$ and the source term $s(w)$ are given by 
\begin{equation*}
f(w) = 
\begin{pmatrix}
\rho \mathbf{u} \\
\rho \mathbf{u} \otimes \mathbf{u} + \frac{1}{M^2} p \mathbb{I}\\
\mathbf{u}(E + p)
\end{pmatrix}
\text{ and }
s(w) = 
\begin{pmatrix}
0\\
\frac{1}{M^2}\frac{\rho}{\alpha}\nabla \beta \\
\frac{\rho}{\alpha} \mathbf{u} \cdot \nabla \beta
\end{pmatrix}.
\end{equation*}
The relaxation model \eqref{sys:full_relax_ELMG} is given by 
\begin{equation}
\label{sys:short_relax}
W_t + \nabla \mathcal{F}(W) = S(W) + \frac{1}{\varepsilon}R(W),
\end{equation}
where $W = (\rho, \rho \mathbf{u},E, \rho \pi, \rho \mathbf{\hat u}, \rho \psi, \rho Z)^T$ denotes the state vector, $\mathcal{F}$ the flux function as defined in \eqref{sys:full_relax_ELMG}. 
The gravitational source term $S(W)$ and the relaxation source term $R(W)$ are given by 
\begin{equation*}
S(W) =
\begin{pmatrix}
0 \\
\frac{1}{M^2} \frac{\rho}{\alpha}\nabla Z \\
\frac{\rho}{\alpha} \mathbf{u} \cdot \nabla Z \\
0 \\
\frac{1}{M^2} \frac{\rho}{\alpha}\nabla Z \\
0 \\
0
\end{pmatrix} 
\text{ and }
R(W) = 
\begin{pmatrix}
0 \\
0 \\
0 \\
\rho(p - \pi) \\
\rho(\mathbf{u - \hat u})\\
\rho(p - \psi) \\
\rho(\beta - Z)
\end{pmatrix}.
\end{equation*}
The relaxation time $\varepsilon$ indicates how fast the perturbed system \eqref{sys:short_relax} is reaching its equilibrium \eqref{sys:short_mod_EulerLMG}.
The relaxation equilibrium state is given by 
\begin{equation}
\label{eq:relaxation_state}
W^{\text{eq}} = \left(\rho, \rho \mathbf{u}, E, \rho p(\rho,e),\rho \mathbf{u},\rho p(\rho,e),\rho \beta\right)^T.
\end{equation}
Following \cite{ChenLevermoreLiu1994}, we can connect \eqref{sys:short_relax} to \eqref{sys:short_mod_EulerLMG} through the matrix $Q \in \mathbb{R}^{(2+d)\times(2(2+d)+1)}$ defined as 
\begin{equation*}
%\label{eq:Q}
Q = 
\begin{pmatrix}
\mathbb{I}_{2+d} & 0_{2(2+d)+1}
\end{pmatrix},
\end{equation*}
where $d$ denotes the dimension.
Then we have for all states $W$ that $QR(W) = 0$ and the physical variables are recovered by $w = QW$ and the flux function $f(w) = Q(\mathcal{F}(W^{\text{eq}}))$.

\section{Time semi-discrete scheme}
To avoid the very restrictive CFL condition that would arise when using an explicit scheme, we will construct an IMEX scheme for which the CFL number is independent of the Mach number.  
Therefore, we split in \eqref{sys:full_relax_ELMG} the flux function $\mathcal{F}(W)$ and source term $S(W)$ in the following way:
\begin{equation}
W_t + \nabla \cdot F(W) + \frac{1}{M^2}\nabla \cdot G(W) = S_E(W) + \frac{1}{M^2}S_I(W) + \frac{1}{\varepsilon} R(W)
\label{eq:decomp}
\end{equation}
where $F(W)$ and $S_E(W)$ will be treated explicitly and $G(W)$ and $S_I(W)$ implicitly. 
The functions $F(W), ~S_E(W), ~G(W)$ and $S_I(W)$ are thus chosen with the purpose of avoiding the need to invert a huge non-linear system which would result treating all terms implicitly.
Instead we propose
\begin{align}
\begin{array}{cc}
F(W) = 
\begin{pmatrix} 
\rho \mathbf{u} \\
\rho \mathbf{u}\otimes\mathbf{u} + \pi \mathbbm{1} + \frac{1-M^2}{M^2} \psi \mathbbm{1} \\
( E + M^2 \pi + (1-M^2)\psi)\mathbf{u} \\
\rho \pi \mathbf{u} + a^2 \mathbf{u}\\
\rho \mathbf{u} \otimes \mathbf{\hat{u}} \\
\rho \psi \mathbf{u}\\
\rho Z \mathbf{u}
\end{pmatrix},&
S_E(W) = 
\begin{pmatrix} 
0 \\
\frac{1}{M^2}\frac{\rho}{\alpha} \nabla Z\\
\frac{\rho}{\alpha} \mathbf{u} \cdot \nabla Z \\
0\\
0\\
0\\
0
\end{pmatrix}, \\
G(W) = 
\begin{pmatrix}
0 \\ 0 \\ 0 \\ 0 \\ \psi\\ a^2 M^2\mathbf{\hat{u}}\\0
\end{pmatrix},
&
S_I(W) = 
\begin{pmatrix}
0 \\ 0 \\ 0 \\ 0 \\ 0\\ \frac{\rho}{\alpha}\nabla Z \\ 0
\end{pmatrix}.
\end{array}
\label{eq:DefSEI}
\end{align} 
The relaxation source term $R$ will merely drive the system to equilibrium, as is standard in relaxation schemes. 
The time semi-discrete scheme is given by the following sequence of implicit, explicit and relaxation steps
\begin{align}
\text{Implicit: }	W_t + \frac{1}{M^2}\nabla \cdot G(W) & = \frac{1}{M^2}S_{I}(W), \label{eq:implicit} \\
\text{Explicit: }~~~~~ 	W_t + \nabla \cdot F(W) & = S_{E}(W), \label{eq:explicit}\\
\text{Projection: }~~~~~~~~~~~~~~~~~~~~~	W_t &= \frac{1}{\varepsilon} R(W). \label{eq:projection}
\end{align}
The projection step \eqref{eq:projection} is equivalent to solving $R(W) = 0$ for $\varepsilon = 0$, see \cite{ChenLevermoreLiu1994}.
Due to the simple structure of $R(W)$, we can immediately set $W = W^{\text{eq}}$ as defined in \eqref{eq:relaxation_state} thus guaranteeing that the data at the new time step is on the equilibrium manifold and thus the original equations \eqref{sys:short_mod_EulerLMG} are satisfied at the new time step.
The formal time semi-discrete scheme is then given by 
\begin{align}
W^{(1)} - W^{n,\text{eq}} + \frac{\Delta t}{M^2} \nabla \cdot G(W^{(1)}) &=\frac{\Delta t}{M^2}S_I(W^{(1)}), \label{eq:t-semi-impl} \\
W^{(2)} - W^{(1)} + \Delta t ~\nabla \cdot F(W^{(1)}) &= \Delta t ~S_E(W^{(1)}), \label{eq:t-semi-expl} \\
W^{n+1} = W^{(2),\text{eq}}\label{eq:t-semi-proj},
\end{align}
Hydrostatic equilibria of  \eqref{eq:t-semi-impl}-\eqref{eq:t-semi-proj} are then given by 
\begin{align}
\begin{rcases}
\mathbf{\hat u}^{(1)} &= 0, \\
\nabla \psi^{(1)} &= \frac{\rho}{\alpha}\nabla Z^{(1)} \\
\end{rcases} \text{Implicit}\label{eq:hydro_impl}\\
\begin{rcases}
\mathbf{u}^{(1)} &= 0, \\
\nabla \pi^{(1)} + \frac{1-M^2}{M^2}\nabla \psi^{(1)} &= \frac{1}{M^2}\frac{\rho}{\alpha}\nabla Z^{(1)}\\
\end{rcases} \text{Explicit}%\notag
\label{eq:hydro_expl}
\end{align}
From \eqref{eq:hydro_impl} and \eqref{eq:hydro_expl} we see that if the implicit step is well-balanced, then the hydrostatic equation for the explicit step reduces to solving
\begin{align*}
\begin{rcases}
\mathbf{u}^{(1)} &= 0, \\
\nabla \pi^{(1)} &= \frac{\rho}{\alpha}\nabla Z^{(1)}\\
\end{rcases}
%	\label{eq:hydro_expl_simpl}
\end{align*}
which is independent of the Mach number. 
%The aim is then to find a space discretization for the implicit step \eqref{eq:t-semi-impl} and explicit step  \eqref{eq:t-semi-expl} that is well-balanced thus fulfilling \eqref{eq:hydro_impl} and \eqref{eq:hydro_expl_simpl} respectively.

\subsection{Mach number expansion of $\psi^{(1)}$}

Due to the sparse structures of the implicit flux function $G$ and implicit source term $S_I$ in \eqref{eq:DefSEI}, the implicit part reduces to solving only two coupled equations in the relaxation variables $\mathbf{\hat{u}},\psi$ given by
\begin{align}
\begin{split}
(\rho\mathbf{\hat{u}})_t + \frac{1}{M^2} \nabla \psi &= \frac{1}{M^2}\kappa\nabla Z, \\
(\rho\psi)_t + a^2 \nabla \cdot \mathbf{\hat{u}} &= 0,
\end{split}
\label{eq:Implicit_part}
\end{align} 
where $\kappa = \rho/\alpha$.
As done in \cite{CordierDegondKumbaro2012,ThomannZenkPuppoKB2019,BispenLukacovaYelash2017}, we rewrite the coupled system \eqref{eq:Implicit_part} into a single equation with an elliptic operator for $\psi$ starting from the time-semi-discrete scheme
\begin{align}
\frac{\rho^{(1)} - \rho^n}{\Delta t} &= 0,\label{eq:time-semi-rho}\\
\frac{(\rho \mathbf{\hat{u}})^{(1)} - (\rho \mathbf{\hat{u}})^{n} }{\Delta t} + \frac{1}{M^2}\nabla \psi^{(1)} - \frac{1}{M^2}\kappa^{(1)}\nabla Z^{(1)}&= 0, \label{eq:time-semi-uhat}\\
\frac{(\rho \psi)^{(1)} - (\rho \psi)^n}{\Delta t} + a^2 \nabla \cdot \mathbf{\hat{u}}^{(1)} &=0. \label{eq:time-semi-psi}
\end{align}
Note, that $Z$, representing the pressure of the steady state, is constant in time and we have $Z^{(1)} = Z^n$, as well as $\alpha^{(1)} = \alpha^n$. 
From the density equation \eqref{eq:time-semi-rho} it follows that $\rho^{(1)} = \rho^n$. 
Together we have $\kappa^{(1)} = \frac{\rho^n}{\alpha^n} = \kappa^n$.
Inserting \eqref{eq:time-semi-uhat} into \eqref{eq:time-semi-psi} we have
\begin{align}
\begin{split}
\psi^{(1)}-\Delta t^2 a^2 \tau^n \nabla \cdot (\tau^n \frac{1}{M^2} \nabla \psi^{(1)}) = \psi^n &- \Delta t^2 a^2 \tau^n \nabla \cdot (\tau^n \frac{\kappa^n}{M^2} \nabla \beta^{n})\\
&- \Delta t a^2 \tau^n \nabla \cdot \mathbf{u}^n,
\end{split}
\label{eq:elliptic_impl}
\end{align}
where we have simplified the notation by using $\tau = 1/\rho$.
Since the data at time $t^n$ is in relaxation equilibrium, we have $\mathbf{\hat{u}}^n = \mathbf{u}^n$ and $Z^n = \beta^n$ on the right hand side of \eqref{eq:elliptic_impl}.
Note that, in contrary to \cite{BispenLukacovaYelash2017,CordierDegondKumbaro2012}, the update \eqref{eq:elliptic_impl} is linear in $\psi$.

Now we analyse the implicit update of $\psi^{(1)}$ with respect to the Mach number.
We assume that the initial data is well-prepared, that is $w^n \in \Omega_{wp}^{\alpha\beta}$ as defined in \eqref{eq:Omegawpab}.
To preserve the scaling of the pressure, we define the following boundary conditions for $\psi$ on the computational domain $D$
\begin{align}
\begin{rcases}
\nabla \psi_0^{(1)} &= \nabla p_0^n\\
\nabla \psi_1^{(1)} &= \nabla p_1^n
\end{rcases}
\text{ on } \partial D.
\label{bc:pressure}
\end{align}
Inserting the Mach number expansion according to $\Omega_{wp}^{\alpha\beta}$ for well-prepared data into \eqref{eq:elliptic_impl} and separating the $\mathcal{O}(M^{-2})$ terms we find
\begin{equation}
\label{eq:implExpOM-2}
\begin{cases}
\nabla \cdot \left(\tau_0^n \nabla \psi_0^{(1)}\right) = \nabla \cdot \left(\tau_0^n \nabla p_0^n\right) &\text{ in } D\\
\nabla \psi_0^{(1)} = \nabla p_0^n &\text{ on } \partial D
\end{cases}.
\end{equation}
This boundary value problem has the unique solution $\nabla \psi^{(1)}_0 = \nabla p_0^n$ on the whole domain $\overline D$.
Substituting the Mach number expansions of $\psi$ and $\tau$ and collecting the $\mathcal{O}(M^{-1})$ terms leads to
\begin{align}
\label{eq:OrderExpInterm}
\begin{split}
&\tau^n_1 \nabla \cdot \left(\tau^n_0 \nabla \psi_0^{(1)}\right) + \tau_0^n \nabla \cdot \left(\tau^n_1 \nabla \psi_0^{(1)} + \tau_0^n \nabla \psi_1^{(1)}\right) \\
= ~&\tau^n_1 \nabla \cdot \left(\tau^n_0 \frac{\rho_0^n}{\alpha} \nabla \beta^n\right) + \tau_0^n \nabla \cdot \left(\tau^n_1 \frac{\rho_0^n}{\alpha} \nabla \beta^n + \tau_0^n \frac{\rho_1^n}{\alpha} \nabla \beta^n\right).
\end{split}
\end{align}
Due to the well-prepared data, we have the relation $\rho_1 \frac{\nabla \beta}{\alpha} = \nabla p_1$ from \eqref{eq:Omegawpab}. 
Then we can simplify the equation \eqref{eq:OrderExpInterm} using $\nabla \psi^{(1)}_0 = \nabla p_0^n$ to
\begin{align}
\label{eq:implExpOM-1}
\begin{cases}
\nabla \cdot \left(\tau_0^n \nabla \psi_1^{(1)}\right) = \nabla \cdot \left(\tau_0^n \nabla p_1\right) &\text{ in } D\\
\nabla \psi^{(1)}_1 = \nabla p_1^n &\text{ on } \partial D
\end{cases}
\end{align}
which has the unique solution $\nabla \psi^{(1)}_1 = \nabla p_1^n$ on the whole domain $\overline D$.
As a last step we look at the $\mathcal{O}(M^0)$ terms and find using the results from \eqref{eq:implExpOM-2} and \eqref{eq:implExpOM-1} that
\begin{equation*}
\nabla \cdot \left( \tau_0^n \nabla \psi^{(1)}_2 \right) =  ~\nabla \cdot \left( \tau_0^n \frac{\rho_2^n}{\alpha} \nabla \beta^n \right) ~~~\text{    in } D.
\label{eq:implExpOM0}
\end{equation*} 
This means the first two terms in the expansion of $\psi^{(1)}$ fulfil the hydrostatic equilibrium \eqref{eq:WPM-2}, \eqref{eq:WPM-1}. 
This proves that the pressure $\psi^{(1)}$ has the correct asymptotic behaviour.

\subsection{Asymptotic preserving property}
Having established the Mach number expansion of $\psi^{(1)}$, we can show now that the time semi-discrete scheme \eqref{eq:t-semi-impl} - \eqref{eq:t-semi-proj} for $M \to 0$ coincides with the time-discretization of the limit equations \eqref{sys:LimitEQ} and that the scheme preserves the set of well-prepared data $\Omega_{wp}^{\alpha\beta}$.
We start by inserting the Mach number expansions given in \eqref{exp:MFr} into \eqref{eq:t-semi-expl}. 
Then we find for the zero order terms in the density, momentum and energy equation as
\begin{align*}
\begin{split}
\rho_0^{n+1} - \rho_0^n + \Delta t ~ \nabla \cdot \rho^n_0 \mathbf{u}_0^n & = 0, \\
\rho_0^{n+1}\mathbf{u}_0^{n+1} - \rho_0^n \mathbf{u}_0^n + \Delta t ~ \left(\rho_0^n \mathbf{u}_0^n \otimes \mathbf{u}^n_0 + \nabla \psi^{(1)}_2 \right) & = \Delta t \frac{\rho_2^n}{\alpha}\nabla \beta^n, \\
\rho^{n+1}_0 e_0^{n+1} - \rho_0^n e_0^n + \Delta t ~\left( \nabla \cdot \mathbf{u}_0^n \left(\rho_0^n e_0^n + \psi^{(1)}_0\right)\right) &= \Delta t ~\frac{\rho_0}{\alpha} \mathbf{u}_0^n \cdot \nabla \beta^n.
\end{split}
\end{align*}
We can simplify the equations by using $\nabla \psi^{(1)}_0 = \nabla p_0^n$ and well-prepared data $w^n \in \Omega_{wp}^{\alpha\beta}$:
\begin{align*}
\begin{split}
\rho_0^{n+1} - \rho_0^n & = 0, \\
\mathbf{u}_0^{n+1} - \mathbf{u}_0^n + \Delta t ~ \left( \mathbf{u}_0^n \cdot \nabla \mathbf{u}_0^n + \frac{\nabla \psi^{(1)}_2}{\rho_0^n} \right) & = \Delta t \frac{\rho_2^n}{\rho_0^n ~\alpha}\nabla \beta^n, \\
p_0^{n+1} - p_0^n &= 0.
\end{split}
\end{align*}
From the first and the last equation we see that $\rho_0$ and $p_0$ do not change in time and looking at the $\mathcal{O}(M^1)$ terms in the energy equation we have $p_1^{n+1} = p_1^n + \mathcal{O}(\Delta t)$. 
This means the pressure and density at $t^{n+1}$ are still well-prepared up to perturbations of $\Delta t$.
Next, we analyse the divergence free property of $\mathbf{u}_0^{n+1}$ and $\rho_0^{n+1}\mathbf{u}_0^{n+1}$. 
This is done by applying the divergence operator on the momentum equation and simplifying using \eqref{eq:implExpOM0}.
We obtain
\begin{align*}
\nabla \cdot \mathbf{u}_0^{n+1} &= \Delta t ~ \nabla \cdot \left( - \mathbf{u}_0^n \cdot \nabla \mathbf{u}_0^n \right) = \mathcal{O}(\Delta t), \\
\nabla \cdot \left(\rho_0^{n+1}\mathbf{u}_0^{n+1}\right) &= \Delta t ~ \nabla \cdot \left(-\rho_0^n \mathbf{u}_0^n \cdot \nabla \mathbf{u}^n_0 - \nabla \psi^{(1)}_2 + \frac{\rho_2^n}{\alpha}\nabla \beta^n\right) = \mathcal{O}(\Delta t).
\end{align*}
For showing the orthogonality condition for $\mathbf{u}_0^{n+1}$ we multiply the momentum equation by $\frac{\nabla \beta}{\alpha}$ and obtain
\begin{equation*}
\mathbf{u}_0^{n+1} \cdot \frac{\nabla \beta^n}{\alpha^n} =  \Delta t ~ \left( - \mathbf{u}_0^n \cdot \nabla \mathbf{u}_0^n - \frac{\nabla \psi^{(1)}_2}{\rho_0^n} + \frac{\rho_2^n}{\rho_0^n ~\alpha}\nabla \beta^n \right) \cdot \frac{\nabla \beta^n}{\alpha} =\mathcal{O}(\Delta t).
\end{equation*}
Therefore all three conditions are satisfied up to a perturbation in $\Delta t$. 
An analogue estimate for the homogeneous case can be found in the method proposed in \cite{CordierDegondKumbaro2012}. 
This analysis yields the following result about the asymptotic preserving property. 
\begin{thm}[AP property]
	For well-prepared initial data $w^n \in \Omega_{wp}^{\alpha\beta}$ and under the boundary conditions \eqref{bc:pressure} the time semi-discrete scheme \eqref{eq:t-semi-impl}- \eqref{eq:t-semi-proj} is asymptotic preserving when $M \to 0$ in the sense that if $w^n \in \Omega_{wp}^{\alpha\beta}$ then also $w^{n+1} \in \Omega_{wp}^{\alpha\beta}$ and in the limit $M \to 0$ the time semi-discrete scheme is a consistent time discretization of the limit equations \eqref{sys:LimitEQ} within $\mathcal{O}(\Delta t)$ terms.
\end{thm}
We remark that the analysis still holds if instead of $\Omega_{wp}^{\alpha\beta}$ the original well-prepared set $\Omega_{wp}$ is used. 

\section{Derivation of the fully discrete scheme}
The derivation of the fully discrete scheme is done in one spatial direction for simplicity, but it can be extended straightforwardly to $d$ dimensions using dimensional splitting in the explicit part and discretizing the expressions 
\begin{equation}
\label{eq:deriv_mulitd}
\begin{split}
\nabla &\cdot (\tau \nabla (\cdot)) = \partial_{x_1}(\tau \partial_{x_1}(\cdot)) + \dots +\partial_{x_d}(\tau \partial_{x_d}(\cdot)) \text{ and } \\
\nabla &\cdot \mathbf{u} = \partial_{x_1} u_1 + \dots +\partial_{x_d} u_d
\end{split}
\end{equation}
with $\mathbf{u} = (u_1, \dots, u_d)$ component-wise in the implicit step.
We use a uniform cartesian grid on a computational domain $D$ divided in $N$ cells $C_i = (x_{i-1/2},x_{i+1/2})$ of step size $\Delta x$. 
We use a standard finite volume setting, where we define at time $t^n$ the piecewise constant functions $w(x,t^n) = w_i^n, \text{ for } x \in C_i.$
\subsection{Well-balanced property of the implicit part}
Applying central differences in \eqref{eq:elliptic_impl} we obtain
\begin{align}
\label{eq:update_impl}
\begin{split}
\psi_i^{(1)} &- \frac{\Delta t^2}{\Delta x^2}\frac{a^2}{M^2} \tau_i^n \left(\tau_{i+1/2}^n (\psi_{i+1}^{(1)} - \psi_i^{(1)}) - \tau_{i-1/2}^n (\psi_{i}^{(1)}-\psi_{i-1}^{n+1})\right) = \\
\psi_i^n &- \frac{\Delta t^2}{\Delta x^2}\frac{a^2}{M^2} \tau_i^n \left(\tau_{i+1/2}^n \kappa_{i+1/2}^n (\beta_{i+1}^{n} - \beta_i^{n}) - \tau_{i-1/2}^n \kappa_{i-1/2}^n(\beta_{i}^{n}-\beta_{i-1}^{n})\right) \\
& - \frac{\Delta t}{2\Delta x} a^2 \left({u}_{i+1}^n - {u}_{i-1}^n\right) ,
\end{split}
\end{align}
where $\tau_{i+1/2} = \frac{1}{2}\left(\tau_{i+1} + \tau_{i}\right)$.
\begin{lem}[Well-balancedness of the implicit part]
	\label{lem:WBpsi}
	Let the initial condition $w^n_i$ be well-balanced, that is
	\begin{equation}
	u_i = 0, ~~~\frac{\rho_i^n}{\alpha_i^n} = 1, ~~~\frac{p_i^n}{\beta_i^n} =1.
	\label{cond:WB}
	\end{equation}
	If the function $\kappa$ is discretized such that in the hydrostatic equilibrium  holds, ie.
	\begin{equation}
	\kappa_{i+1/2} = 1,
	\label{cond:K}
	\end{equation}
	then it is $\psi^{(1)}_i = \psi^n_i$ for all cells $i = 1, \dots N$, that means \eqref{eq:t-semi-impl} is well-balanced in the sense that $W^{(1)}$ fulfils \eqref{eq:hydro_impl}.
\end{lem}
\begin{pf}
	From the condition \eqref{cond:WB} we have $\kappa_{i+1/2} = 1$.
	At time level $t^n$ we know that $\psi^n = p^n$.
	Therefore we can write
	\begin{equation}
	\psi_{i+1}^n - \psi_i^n = \beta_{i+1}^n - \beta_i^n = \kappa_{i+1/2}^n (\beta_{i+1}^n - \beta_i^n).
	\label{eq:psinWB}
	\end{equation}
	Using $u = 0$ and inserting \eqref{eq:psinWB} into \eqref{eq:update_impl}, we have 
	\begin{align}
	\label{eq:impl_upd_wb}
	\begin{split}
	\psi_i^{(1)} &- \frac{\Delta t^2}{\Delta x^2}\frac{a^2}{M^2} \tau_i^n \left(\tau_{i+1/2}^n (\psi_{i+1}^{(1)} - \psi_i^{(1)}) - \tau_{i-1/2}^n (\psi_{i}^{(1)}-\psi_{i-1}^{n+1})\right) = \\
	\psi_i^n &- \frac{\Delta t^2}{\Delta x^2}\frac{a^2}{M^2} \tau_i^n \left(\tau_{i+1/2}^n (\psi_{i+1}^{n} - \psi_i^{n}) - \tau_{i-1/2}^n (\psi_{i}^{n}-\psi_{i-1}^{n})\right).
	\end{split}
	\end{align}
	Define the tridiagonal coefficient matrix $A$ by 
	\begin{equation*}
	A = \mathrm{diag}(-\mu\tau_i^n \tau_{i-1/2}^n, 1 +\mu \tau_i^n(\tau_{i+1/2}^n+\tau_{i-1/2}^n),-\mu\tau_i^n \tau_{i+1/2}^n),
	\end{equation*} 
	where $\mu = \frac{\Delta t^2}{\Delta x^2}\frac{a^2}{M^2}$.
	Then we can write \eqref{eq:impl_upd_wb} as
	\begin{align}
	A \psi^{(1)} = A \psi^n \Leftrightarrow A (\psi^{(1)} - \psi^n) = 0,
	\label{eq:psiWB}
	\end{align}
	Since the matrix $A$ is strict diagonal dominant it is invertible. 
	Then we have from \eqref{eq:psiWB} that $\psi_i^{(1)} = \psi_i^n$ for all $i = 1, \dots, N$.
	The proof can be extended to $d$ dimensions using \eqref{eq:deriv_mulitd} for the space discretization.
	In $d$ dimensions the coefficient matrix $A$ is an invertible strict diagonal dominant banded Matrix with $2d+1$ diagonals. 
	Therefore the results holds also in $d$ dimensions. 
\end{pf}

In the following we will use a second order accurate discretization of $\kappa_{i+1/2}$ that fulfils \eqref{cond:K} and is given by
\begin{equation*}
\kappa_{i+1/2}= \frac{1}{2}\left(\frac{\rho_{i+1}}{\alpha_{i+1}} + \frac{\rho_i}{\alpha_i}\right).
\end{equation*}

\subsection{Godunov type finite volume scheme} 
We consider the explicit step \eqref{eq:t-semi-expl} using the explicit operators $F$ and $S_E$ defined in \eqref{eq:DefSEI}.
\begin{align}
\begin{split}
\partial_t \rho + \partial_{x} \rho u &= 0 \\
\partial_t \rho u + \partial_{x} (\rho u^2 + \pi + \frac{1-M^2}{M^2}\psi) &= \frac{1}{M^2}\kappa\partial_{x}Z\\
\partial_t E + \partial_{x}((E + M^2 \pi + (1-M^2)\psi) u) &= u \kappa\partial_{x}Z \\
\partial_t \rho \pi + \partial_{x}((\rho\pi + a^2)u) &= 0 \\
\partial_t \rho\hat{u} + \partial_{x}(\rho u \hat{u}) &= 0 \\
\partial_t \rho\psi + \partial_{x}(\rho \psi u) &= 0\\
\partial_t \rho Z + \partial_{x}(\rho Z u) &=0.
\end{split}
\label{sys:expl}
\end{align}
The derivation of the Godunove type finite volume scheme follows closely the steps given eg. in \cite{ThomannZenkPuppoKB2019,BerthonKlingenbergZenk2018,DesveauxZenkBerthonKlingenberg2016,ThomannZenkKlingenberg2019,Toro2009}. 
The omitted proofs to the results given in this section can be done analogously following those references. 
To construct a Riemann solver for \eqref{sys:expl}, we follow \cite{DesveauxZenkBerthonKlingenberg2016} and include the source term in the flux formulation. 
To calculate the Riemann invariants given in Lemma \ref{lem:EW}, we rewrite \eqref{sys:expl} in non-conservative form using the primitive variables $(\rho, u, e, \pi, \hat u, \psi, Z)$.
Since Riemann invariants are invariant under change of variables, they are the same as for the equations in conservation form. 

\begin{lem}
	\label{lem:EW}
	System \eqref{sys:expl} admits the linear degenerate eigenvalues $\lambda^\pm =u \pm \frac{a}{\rho}$ and $\lambda^u = u$, where the eigenvalue $\lambda^u$ has multiplicity 5.
	The relaxation parameter $a$ as well as the eigenvalues are independent of the Mach number $M$. 
	The Riemann invariants with respect to $\lambda^u$ are
	\begin{equation*}
	%	\label{eq:RI_u}
	I_1^u = u,~ I_2^u = M^2 \pi + (1 - M^2) \psi - \kappa Z
	\end{equation*}
	and with respect to $\lambda^\pm$ 
	\begin{align*}
	%	\label{eq:RI_pm}
	\begin{split}
	I_1^\pm &= u \pm \frac{a}{\rho},~I_2^\pm = \pi + \frac{a^2}{\rho}, \\
	I^\pm_3 &= e - \frac{M^2}{2a^2}\pi^2 - \frac{1-M^2}{a^2} \pi \psi,\\ 
	I_{4}^\pm &= \hat{u},~I^\pm_5 = \psi, ~I^\pm_6 = Z.
	\end{split}
	\end{align*}
\end{lem}
We will follow the theory of Harten, Lax and van Leer \cite{HartenLaxVanLeer1983} for deriving an approximate Riemann solver $W_\mathcal{RS}\left(\frac{x}{t}; W_L^{(1)}, W_R^{(1)}\right)$ based on the states $W^{(1)}$ after the implicit step.
Due to the linear-degeneracy from Lemma \ref{lem:EW}, the structure of the approximate Riemann solver is given as follows
%, as displayed in Figure \ref{fig:Riemann_solution},
\begin{align}
W_\mathcal{RS}\left(\frac{x}{t}; W_L^{(1)}, W_R^{(1)}\right) = 
\begin{cases}
W_L^{(1)} & \frac{x}{t} < \lambda^- ,\\
W^*_L & \lambda^- < \frac{x}{t} < \lambda^u ,\\
W^*_R & \lambda^u < \frac{x}{t} < \lambda^+ ,\\
W_R^{(1)} & \lambda^+ < \frac{x}{t}.
\end{cases} 
\label{eq:Sol_RP}
\end{align}
%\begin{figure}[htbp]
%	\centering
%	\begin{tikzpicture}[scale=2.]
%	\begin{scope}[every node/.style={scale=1.}]
%	\draw
%	(-0.5,0)-- (2.5,0);
%	\draw (0,0.85) coordinate (l1) node[above] {$u - a/\rho$} -- (1,0) coordinate (x0) node[below, yshift=-0.5ex] {$x=0$} ;
%	\draw (2,0.85) coordinate (l3) node[above right] {$u + a/\rho$} -- (1,0);
%	\draw (1.15,0.85) coordinate (l2) node[above,yshift=0.5ex] {$u$} -- (1,0);
%	\draw (0.25,0.4) coordinate (WL) node[left] {$W_L$};
%	\draw (1.75,0.4) coordinate (WR) node[right] {$W_R$};
%	\draw (0.7,0.5) coordinate (WSL) node[above] {$W^*_L$};
%	\draw (1.4,0.5) coordinate (WSR) node[above] {$W^*_R$};
%	\end{scope}
%	\end{tikzpicture}
%	\caption{Structure of the Riemann solution.}
%	\label{fig:Riemann_solution}
%\end{figure}
To compute the intermediate states $W_{L,R}^{\ast}$, we use the Riemann invariants as given in Lemma \ref{lem:EW}. 

%Note that since the eigenvalues $\lambda^\pm$ have multiplicity 1, we get the expected 6 Riemann invariants.
%This does not hold in general for eigenvalues with multiplicity larger than 1. 
%Nevertheless, the invariants \eqref{eq:RI_u} and \eqref{eq:RI_pm} give enough relations to determine the solution to a Riemann problem for \eqref{sys:expl} as shown in the following lemma. 
\begin{lem}
	\label{lem:RSolver}
	Consider an initial value problem with initial data $W = W^{(1)}$ given by 
	\begin{align*}
	W_0(x)=\begin{cases}
	W_L &x < 0\\
	W_R &x > 0
	\end{cases}.
	\end{align*}
	Then, the solution consists of four constant states separated by contact discontinuities with the structure given in \eqref{eq:Sol_RP}. 
	The solution for the intermediate states $W^\ast_{L},W^\ast_R$ with $u^\ast= u_L^\ast = u_R^\ast$ is given by 
	\begin{align}
	\begin{split}
	\frac{1}{\rho_L^\ast} &= \frac{1}{\rho_L} + \frac{1}{a^2}(\pi_L - \pi_L^\ast),\\ \frac{1}{\rho_R^\ast} &= \frac{1}{\rho_R} + \frac{1}{a^2}(\pi_R - \pi_R^\ast),\\
	u^\ast & =  \frac{1}{2}(u_{L} + u_{R}) - \frac{1}{2a}\left(\pi_R -\pi_L + \frac{1-M^2}{M^2}(\psi_R - \psi_L) -\frac{\kappa}{M^2}(Z_R - Z_L)\right),\\
	\pi_{L}^\ast &= \frac{1}{2}(\pi_L + \pi_R) - \frac{a}{2}(u_{R} - u_{L}) + \frac{1-M^2}{2M^2}(\psi_R - \psi_L) - \frac{\kappa}{2 M^2}(Z_R - Z_L), \\
	\pi_{R}^\ast &= \frac{1}{2}(\pi_L + \pi_R) - \frac{a}{2}(u_{R} - u_{L}) - \frac{1-M^2}{2M^2}(\psi_R - \psi_L) + \frac{\kappa}{2 M^2}(Z_R - Z_L), \\
	e_{L}^\ast &= e_{L} - \frac{1}{2 a^2}(\pi_{L}^2 - (\pi_{L}^\ast)^2 + (1-M^2)(\pi_{L} - \pi_{L}^\ast)\psi_{L}), \\
	e_{R}^\ast &= e_{R} - \frac{1}{2 a^2}(\pi_{R}^2 - (\pi_{R}^\ast)^2 + (1-M^2)(\pi_{R} - \pi_{R}^\ast)\psi_{R}), \\
	\psi_{L,R}^\ast &=\psi_{L,R},\\
	\hat{u}_{L,R}^\ast &= \hat{u}_{L,R},\\
	Z_{L,R}^\ast &= Z_{L,R}. 
	\end{split}
	\label{eq:intermStates}
	\end{align}
\end{lem}
Having established the structure of the Riemann solver, we can show that it is preserving hydrostatic equilibria.
\begin{lem}[Well-balancedness of Riemann Solver]
	\label{lem:WBRS}
	Let the initial condition $w^n_L, w^n_R$ be given in hydrostatic equilibrium \eqref{cond:WB}. 
	Let the function $\kappa$ be defined as in \eqref{cond:K}. 
	Then the intermediate states \eqref{eq:intermStates} satisfy 
	\begin{equation*}
	W_L^{(1)\ast} = W_L^{(1)}, ~W_R^{(1)\ast} = W_R^{(1)}
	\end{equation*}
	that is, the approximate Riemann solver as defined in Lemma \ref{lem:RSolver} is at rest.
\end{lem}
\begin{pf}
	From Lemma \ref{lem:WBpsi}, we know that $\psi^{(1)} = p^n$ and satisfies
	\begin{equation}
	\psi_L^{(1)} - \psi_R^{(1)} = \kappa(Z_L^{(1)} - Z_R^{(1)}).
	\label{eq:lemWB_psi}
	\end{equation} 
	We also know that $\pi^{(1)}_{L,R} = \pi^n_{L,R} = p^n_{L,R}$ and since $w^n$ is fulfilling \eqref{cond:WB} and with \eqref{cond:K} we have $\pi^{(1)}_L - \pi^{(1)}_R = \kappa(Z_L^{(1)} - Z_R^{(1)})$.
	Then we have 
	\begin{align*}
	&\pi_R^{(1)} - \pi_L^{(1)} + \frac{1-M^2}{M^2}(\psi_R^{(1)} - \psi_L^{(1)}) - \frac{\kappa}{M^2}(Z_R^{(1)} - Z_L^{(1)}) = \\
	&\pi_R^{(1)} - \pi_L^{(1)} - \kappa(Z_R^{(1)} - Z_L^{(1)}) = 0.
	\end{align*}
	Since $u^{(1)}_{L,R} = u^n_{L,R} = 0$, we find $u^{(1)\ast} = 0$. 
	With $u^{(1)\ast} = 0$ and \eqref{eq:lemWB_psi} and the fact that $\psi^n = p^n = \pi^n = \pi^{(1)}$, we can write
	\begin{align*}
	\pi_{L}^\ast &= \frac{1}{2}(\pi_L^{(1)} + \pi_R^{(1)}) + \frac{1-M^2}{2M^2}(\psi_R^{(1)} - \psi_L^{(1)}) - \frac{\kappa}{2 M^2}(Z_R^{(1)} - Z_L^{(1)})\\
	&= \frac{1}{2}(\pi_L^{(1)} + \pi_R^{(1)}) - \frac{1}{2}(\pi_R^{(1)} - \pi_L^{(1)})\\
	&= \pi_L.
	\end{align*}
	Analogously follows $\pi_R^\ast = \pi_R$.
	Then it follows directly from the intermediate states \eqref{eq:intermStates} that $\rho_L^\ast = \rho_L$, $\rho_R^\ast = \rho_R$ and $e_L^\ast = e_L$, $e_R^\ast = e_R$.
\end{pf}
\vskip0.25cm
Another important property is that the density and pressure remain positive during the simulation. 
This is equivalent to preserving the following domain
\begin{equation*}
\Omega_{phy} = \left\lbrace w \in \Omega, \rho > 0, e > 0 \right\rbrace.
\end{equation*}
We show that the Riemann solver preserves $\Omega_{phy}$.
\begin{lem}[Positivity preserving property of Riemann Solver]
	\label{lem:posRS}
	Suppose the initial data  $W^{(1)}_{L,R}$ is composed of $w_{L,R}^{(1)} \in \Omega_{phy} \cup ~ \Omega_{wp}^{\alpha,\beta}$ and $\psi^{(1)}$ satisfies the boundary conditions \eqref{bc:pressure}. 
	Then solution of the Riemann problem given by $Q W_{\mathcal{RS}}(\frac{x}{t};W_L^{(1)},W_R^{(1)})$ is contained in $\Omega_{phy}$ for a relaxation parameter $a$ sufficiently large but independent of $M$.
\end{lem}
\begin{pf}
	The proof for the intermediate states for the density can be taken from \cite{ThomannZenkPuppoKB2019,ThomannZenkKlingenberg2019}. 
	After the implicit step we have $u^{(1)} = u^n$, $\pi^{(1)} = \pi^n$ and $Z^{(1)} = Z^n$.
	We use the following notation $\Delta (\cdot) = (\cdot)_R - (\cdot)_L$.
	For the internal energy, the intermediate state $\pi_L^{(1)\ast}$ is inserted into $e_L^\ast$ and we have
	\begin{align}
	\label{eq:internal_e_pos}
	\begin{split}
	e_L^{(1)\ast} = &~e_L^n + \frac{1}{8}\Delta u^2\\
	&+ \frac{1}{2 a^2} \left(-~\left(\pi_L^n\right)^2 + \frac{1}{4}\left(\pi_L^n + \pi_R^n -\Delta \psi^{(1)}+ \frac{1}{M^2}H^{(1)}\right)^2 \right.\\
	& \hspace{1.5cm}\left.+ ~\frac{1}{2}~\psi_L^{(1)}(1-M^2) \left(\Delta \pi^n - \Delta \psi^{(1)}+ \frac{1}{M^2}H^{(1)}\right)\right) \\
	&+\frac{1}{4a}\Delta u^n\left(\Delta \pi^n + 2\pi_L^n - \Delta \psi^{(1)} + \frac{1}{M^2}H^{(1)} + (1-M^2)\psi_L^{(1)}\right),
	\end{split}
	\end{align}
	where we have defined $H^{(1)} = (\psi_R^{(1)} - \psi_L^{(1)}) - \kappa(Z_R^n - Z_L^n)$.
	We know from the Mach number analysis in Section 4 that $\psi^{(1)}$ preserves the hydrostatic equilibrium up to a perturbation of $M^2$, thus $H^{(1)} = \mathcal{O}(M^2)$.
	Therefore we find a relaxation parameter $a > \rho \sqrt{\partial_\rho p(\rho,e)}$ independent of $M$ that can control the negative terms in \eqref{eq:internal_e_pos} and we have $e_L^{(1)\ast} > 0$.
\end{pf}

With the solution of the Riemann problem \eqref{eq:Sol_RP} we can define the numerical fluxes at the interface $x_{i+1/2}$. 
With $S_{i+1/2} = (0,s_{i+1/2},u_i^\ast s_{i+1/2})$ where $s_{i+1/2} = \kappa_{i+1/2}(Z_{i+1} - Z_i)$ we have
\begin{align}
\begin{split}
F^-_{i+1/2} = 
\begin{cases}
F\left(W_i^{(1)}\right), &\lambda^- > 0 \\
F\left(W_i^{(1)\ast}\right), & \lambda^u > 0 > \lambda^- \\
F\left(W_i^{(1)\ast}\right), & \lambda^u  = 0 \\
F\left(W_{i+1}^{(1)\ast}\right)-S_{i+1/2}, & \lambda^+ > 0 > \lambda^u \\
F\left(W_{i+1}^{(1)}\right)-S_{i+1/2}, & \lambda^+ < 0
\end{cases},\\
F^+_{i+1/2} = 
\begin{cases}
F\left(W_i^{(1)}\right)+S_{i+1/2}, &\lambda^- > 0 \\
F\left(W_i^{(1)\ast}\right) + S_{i+1/2}, & \lambda^u > 0 > \lambda^- \\
F\left(W_{i+1}^{(1)\ast}\right), & \lambda^u  = 0 \\
F\left(W_{i+1}^{(1)\ast}\right), & \lambda^+ > 0 > \lambda^u \\
F\left(W_{i+1}^{(1)}\right), & \lambda^+ < 0
\end{cases},
\end{split}
\label{def:F_interface}
\end{align}
where the superscript $(1)$ emphasizes that the states after the implicit step are used.
We want to stress that we include the source term into the flux definition and therefore in general it is $F^-_{i+1/2}\neq F_{i+1/2}^+$.
This leads to the following update of the explicit part
\begin{equation}
W^{(2)}_i = W^{(1)}_i - \frac{\Delta t}{\Delta x} (F^-_{i+1/2} - F^+_{i-1/2}).
\label{eq:update_explicit}
\end{equation}
To avoid interactions between the approximate Riemann solvers at the interfaces $x_{i+1/2}$, we have a CFL restriction on the time step of 
\begin{equation}
\Delta t \leq \frac{1}{2} \frac{\Delta x}{ \underset{i}{\max}|u_i \pm a/\rho_i|}
\label{eq:cfl}
\end{equation} 
which is independent of the Mach number.
Due to the relaxation step \eqref{eq:t-semi-proj}, we can directly give the update of the physical variables $w$ as 
\begin{align}
\label{rem:update_physVar}
\begin{split}
w_i^{n+1} = w_i^n - \frac{\Delta t}{\Delta x}&\left(QF_{i+1/2}^-\left(W_{\mathcal{RS}}\left(0;W_i^{(1)},W_{i+1}^{(1)}\right)\right) \right.\\
&~~~\left.-~ QF_{i-1/2}^+\left(W_{\mathcal{RS}}\left(0;W_{i-1}^{(1)},W_{i}^{(1)}\right)\right)\right).
\end{split}
\end{align}
\begin{thm}[Well-balanced property 1]
	Let $w_i$ on all cells $i \in \left\lbrace1,N\right\rbrace$ be given in hydrostatic equilibrium \eqref{cond:WB}. Let $\kappa$ be defined as in \eqref{cond:K}.
	Then the first order scheme given by the steps \eqref{eq:update_impl},\eqref{rem:update_physVar} is well-balanced.
\end{thm}
\begin{pf}
	Since $w^n$ fulfils the hydrostatic equilibrium, we know from Lemma \ref{lem:WBpsi} that $W^{(1)}_i = W^n_i$ fulfils the hydrostatic equilibrium and from Lemma \ref{lem:WBRS} that the approximate Riemann solver at the cell interfaces is at rest. 
	With the definition of the fluxes \eqref{def:F_interface}, we have 
	\begin{equation*}
	F_{i-1/2}^+ = F(W_i^n), \hskip 1cm F_{i+1/2}^- = F(W_i^n).
	\end{equation*}
	Using the formulation \eqref{rem:update_physVar} for the update of the variables $w$, we have
	\begin{equation*}
	\label{Theo:WB_update1D}
	w_i^{n+1} = w_i^n - \frac{\Delta t}{\Delta x} Q\left(F_{i+1/2}^- - F_{i-1/2}^+\right) = w_i^n.
	\end{equation*}
	This shows the well-balanced property in one dimension. 
	Since we apply dimensional splitting in the multi-dimensional set-up, the proof can be easily extended by giving the update \eqref{Theo:WB_update1D} as a sum of the flux differences along each dimension.
\end{pf}

\begin{thm}[Positivity preserving 1]
	\label{Theo:Pos1}
	Let the initial state in $d$ dimensions be given as 
	\begin{equation*}
	w_i^n \in \Omega = \Omega_{phy} \cap \Omega_{wp}^{\alpha\beta}
	\end{equation*}
	Then under the Mach number independent CFL condition 
	\begin{equation*}
	\frac{\Delta t}{\Delta x} \underset{i}{\max}|\lambda^\pm(w_i^n)| < \frac{1}{2^d},
	\end{equation*}
	and the boundary conditions \eqref{bc:pressure} the numerical scheme defined by \eqref{eq:update_impl},\eqref{rem:update_physVar} preserves the positivity of density and internal energy, that is $w_i^{n+1} \in \Omega_{phy}$ for a sufficiently large relaxation parameter $a$ independent of $M$.
\end{thm}

An important property for any low Mach scheme is the behaviour of the diffusion. 
Due to the fact that $\psi^{(1)}$ is still well-prepared after the implicit step, the diffusion of the scheme is of order $\mathcal{O}(M^0)$. 
The computations are performed analogously to the homogeneous case and can be found in \cite{ThomannZenkPuppoKB2019}.

\subsection{Second order extension}
Here, we give a strategy to extend the first order scheme to second order accuracy such that the well-balanced and the positivity preserving property are maintained. 

For the time integration, we use the second order scheme presented in \cite{ThomannZenkPuppoKB2019}. 
The second order extension in space is realized by a linear reconstruction of the interface values. 
We reconstruct in the primitive variables $w^p = (\rho,\mathbf{u},p)$ and $\psi$ on each cell.
Since we use dimensional splitting, we reconstruct along each space dimension separately.
We consider a linear function on $C_i$ defined as
\begin{equation}
\label{eq:lin_recon}
w^p(x) = w_i^p + \sigma(x - x_i).
\end{equation}
The slopes $\sigma = (\sigma^\rho, \sigma^u, \sigma^p)$ are obtained by using information from the neighbouring cells. 
The interface values on cell $C_i$ denoted by $w_{i-1/2}^+, w_{i+1/2}^-$ are then obtained by evaluating $w^p(x)$ at the cell interfaces. 
The reconstruction \eqref{eq:lin_recon} has to fulfil two properties.
Firstly, the interface values in conserved variables have to be in $\Omega_{phy}$ to satisfy the conditions in Lemma \ref{lem:posRS}. 
Secondly, if $w^n$ fulfils the hydrostatic equilibrium, also the interface values have to fulfil the hydrostatic equilibrium. 
To meet the first requirement we apply on the slopes $\sigma$ a limiting procedure described in \cite{ThomannZenkKlingenberg2019} to guarantee $w_{i+1/2}^-, w_{i-1/2}^+ \in \Omega_{phy}$.
For the well-balanced property, we apply a hydrostatic reconstruction on the pressure as it can be found in \cite{KaeppeliMishra2016,ThomannZenkKlingenberg2019}
\begin{align}
\label{eq:Ptrafo}
\begin{split}
q_{i-1} &= \pi_{i-1} + s_{i-1/2}, \\
q_{i+1} &= \pi_{i+1} - s_{i+1/2}.
\end{split}
\end{align}
The slope for $\pi$ is then calculated as
\begin{equation*}
\sigma^q = \text{minmod}\left(\frac{q_{i+1} - \pi_i}{\Delta x}, \frac{\pi_i - q_{i-1}}{\Delta x}\right).
\end{equation*}
Analogously we get the modified slope for $\psi^{(1)}$.
This results into $\pi_{i+1/2}^- = \pi_{i-1/2}^+ = p_i^n$ and $\psi_{i+1/2}^{(1),-} = \psi_{i-1/2}^{(1),+} = p_i^n$ when being in a hydrostatic equilibrium and the Riemann Solver is at rest. 
We will summarize the well-balanced and positivity preserving property of the second order scheme. The proofs are analogous to the ones shown in \cite{ThomannZenkPuppoKB2019,ThomannZenkKlingenberg2019}.
\begin{thm}[Well-balanced property 2]
	Let the initial condition $w^n$ be given in hydrostatic equilibrium \eqref{cond:WB}. Let the function $\kappa$ be defined as in \eqref{cond:K}. Then, using the transformation \eqref{eq:Ptrafo}, the second order scheme is well-balanced.
\end{thm}

\begin{thm}[Positivity property 2]
	Let the initial state be given as $w_i^n \in \Omega$ satisfying the boundary conditions \eqref{bc:pressure} and the limiting procedure given in \cite{ThomannZenkKlingenberg2019} is used. Then for a sufiiciently large relaxation parameter $a$, under the Mach number independent CFL condition 
	\begin{equation*}
	\frac{\Delta t}{\Delta x} \underset{i}{\max}|\lambda^\pm(w_i^n)| < \frac{1}{2 \cdot 2^d},
	\end{equation*}
	where $d$ denotes the dimension, the second order scheme preserves the domain  $\Omega_{phy}$.
\end{thm}

\section{Numerical results}
In this section, we give numerical test cases to validate the theoretical properties of the first and second order scheme. 
For all test cases we assume an ideal gas law $p = (\gamma - 1) \rho e$. 
The implicit non-symmetric linear system given by \eqref{eq:update_impl} is solved with the GMRES algorithm combined with a preconditioner based on an incomplete LU decomposition. 
To choose the relaxation parameter $a$, we follow the procedure given in \cite{Bouchut2004} to obtain a local estimate 
for $a$.
We calculate a global estimate by taking the maximum of the local values of $a$ and multiply by a constant $c_a$ independent of $M$ to ensure the stability property given in Lemma \ref{lem:properties}.

\subsection{Well-balanced test case}
%\color{green} Markus suggest to run them longer \color{black} 
To numerically verify the well-balanced property of the scheme, we compute an isothermal equilibrium with a linear potential in two dimensions as given in \eqref{eq:isoth} where $\mathbf{u} = (u_1, u_2), ~\chi = 1$ and $\gamma = 1.4$. 
In Table \ref{tab:WB2Diso} we give the error at the final time $T_f = 1$ for different Mach and Froude numbers  on the domain $D=[0,1]^2$.
The results are computed with the first order scheme.
As expected, the error is of order of machine precision as can be seen in Table \ref{tab:WB2Diso}. 
\begin{table}[htpb]
	\centering
	\renewcommand{\arraystretch}{1.25}
	\begin{tabular}{c c cccc}
		$M$&$Fr$&$\rho$&$\rho u_1$&$\rho u_2$&$E$\\
		\hline\hline 
		$10^{-1}$&$10^{-1}$& 2.459E-017 & 3.605E-016  & 3.605E-016 &  2.419E-017\\
		$10^{-2}$&$10^{-2}$&  5.606E-017 &  9.999E-017 &  9.999E-017 & 5.507E-017\\
		$10^{-3}$&$10^{-3}$& 2.506E-017 &  9.811E-016 &  9.811E-016 &  2.457E-017\\
		$10^{-4}$&$10^{-4}$& 2.539E-017 &  5.304E-017 &  5.304E-017 &  2.495E-017\\
		%		$10^{-4}$&$10^{-1}$& 2.4413804311507194E-017   3.4055459316974665E-016   3.4055459316974665E-016   2.4746871218894741E-017\\
		%		$10^{-1}$&$10^{-4}$&\\
		\hline
	\end{tabular}
	\caption{$L^1$-error of isothermal equilibrium at $T=1$ (non-dimensional).}
	\label{tab:WB2Diso}
\end{table}
\subsection{Accuracy}
To numerically validate the second order accuracy of the proposed scheme, we compare the numerical solution obtained with the second order scheme to an exact solution of the Euler equations with gravity as given in \cite{ChandrashekarKlingenberg2015}.
In physical variables, it is given in 2 dimensions with $\mathbf{x} = (x_1,x_2)$ and $\mathbf{u} = (u_1,u_2)$ as
\begin{align}
\label{eq:ExactInitDim}
\begin{split}
\rho(\mathbf{x},t) &= 1 + 0.2 \sin\left(\pi (x_1 + x_2 - t(u_{1_0} + u_{2_0}))\right) \frac{kg}{m^3} \\
u_1(\mathbf{x},t) &= u_{1_0} \frac{m}{s} \\
u_2(\mathbf{x},t) &= u_{2_0} \frac{m}{s} \\
p(\mathbf{x},t) &= p_0 + t(u_{1_0} + u_{2_0}) - (x_1 + x_2) \\
&\quad\quad~+ 0.2 \cos\left(\pi(x_1 + x_2 - t(u_{1_0} + u_{2_0}))\right)/\pi \frac{kg}{m s^2}.
\end{split}
\end{align} 
For the parameters we set $u_{1_0} = 20, u_{2_0} = 20$ and $p_0 = 4.5$.
The gravitational potential is linear and given as $\Phi(\mathbf{x}) = x_1 + x_2$.
For $\mathbf{u} = 0$, \eqref{eq:ExactInitDim} is in hydrostatic equilibrium and we set $\alpha$ and $\beta$ as the density and pressure of the stationary state respectively.
We want to remark that this equilibrium is neither isothermal nor polytropic.
The computational domain is $D = [0,1]^2$ and the final time $T = 0.01 s$.

\begin{figure}[htpb]
	\centering
	\begin{subfigure}[t]{0.5\textwidth}
		\begin{tikzpicture}[scale=0.8]
		%loglogaxis
		\begin{loglogaxis} [legend style={
			at={(1.4,1)},
			anchor=north}, 
		legend columns = 1,
		title={Density $\left[\frac{kg}{m^3}\right]$},
		xlabel={N},
		ylabel={$L^1$ error},
		xmin = 20, xmax = 250,
		xtick={25,50,100,200},
		xticklabels={$25$,$50$,$100$,$200$}
		]
		\printslope{2}{50}{1E-4}{30}
		% 0.1 0.1
		\addplot table [x=N, y=r, col sep=space] {M0p10p1.dat}; 
		\addplot table [x=N, y=r, col sep=space] {M0p010p01.dat}; 
		\addplot table [x=N, y=r, col sep=space] {M0p0010p001.dat}; 
		\addplot table [x=N, y=r, col sep=space] {M0p00010p0001.dat}; 
		\addplot table [x=N, y=r, col sep=space] {M0p00010p1.dat}; 
		\addplot table [x=N, y=r, col sep=space] {M0p10p0001.dat};
		\end{loglogaxis} 
		\end{tikzpicture}
	\end{subfigure}
	\begin{subfigure}[t]{0.45\textwidth}
		\begin{tikzpicture}[scale=0.8]
		%loglogaxis
		\begin{loglogaxis} [legend style={
			at={(1.4,1)},
			anchor=north}, 
		legend columns = 1,
		title={Momentum $\left[\frac{kg}{m^2 s}\right]$},
		xlabel={N},
		ylabel={$L^1$ error},
		xmin = 20, xmax = 250,
		xtick={25,50,100,200},
		xticklabels={$25$,$50$,$100$,$200$}
		]
		\printslope{2}{50}{2E-3}{30}
		% 0.1 0.1
		\addplot table [x=N, y=m1, col sep=space] {M0p10p1.dat}; 
		\addplot table [x=N, y=m1, col sep=space] {M0p010p01.dat}; 
		\addplot table [x=N, y=m1, col sep=space] {M0p0010p001.dat}; 
		\addplot table [x=N, y=m1, col sep=space] {M0p00010p0001.dat}; 
		\addplot table [x=N, y=m1, col sep=space] {M0p00010p1.dat}; 
		\addplot table [x=N, y=m1, col sep=space] {M0p10p0001.dat}; 
		\end{loglogaxis} 
		\end{tikzpicture}
	\end{subfigure}
	\flushleft\hskip0.4cm
	\begin{subfigure}[t]{0.5\textwidth}
		\begin{tikzpicture}[scale=0.8]
		%loglogaxis
		\begin{loglogaxis} [legend style={
			at={(1.75,0.85)},
			anchor=north}, legend columns = 1, legend style={font=\Large},
		title={Energy $\left[\frac{kg}{m s^2}\right]$},
		xlabel={N},
		ylabel={$L^1$ error},
		xmin = 20, xmax = 250,
		xtick={25,50,100,200},
		xticklabels={$25$,$50$,$100$,$200$}
		]
		\printslope{2}{50}{4E-2}{30}
		% 0.1 0.1 
		\addplot table [x=N, y=E, col sep=space] {M0p10p1.dat}; 
		\addplot table [x=N, y=E, col sep=space] {M0p010p01.dat}; 
		\addplot table [x=N, y=E, col sep=space] {M0p0010p001.dat}; 
		\addplot table [x=N, y=E, col sep=space] {M0p00010p0001.dat}; 
		\addplot table [x=N, y=E, col sep=space] {M0p00010p1.dat}; 
		\addplot table [x=N, y=E, col sep=space] {M0p10p0001.dat}; 
		\legend{$M = 10^{-1} ~Fr=10^{-1}$, $M=10^{-2} ~Fr=10^{-2}$, $M=10^{-3} ~Fr=10^{-3}$, $M=10^{-4} ~Fr=10^{-4}$, $M=10^{-1} ~Fr=10^{-4}$, $M=10^{-4} ~Fr=10^{-1}$}
		\end{loglogaxis} 
		\end{tikzpicture}
	\end{subfigure}
	\caption{$L^1$ error curves in dependence of Mach and Froude number (dimensional).}
	\label{fig:error}
\end{figure}
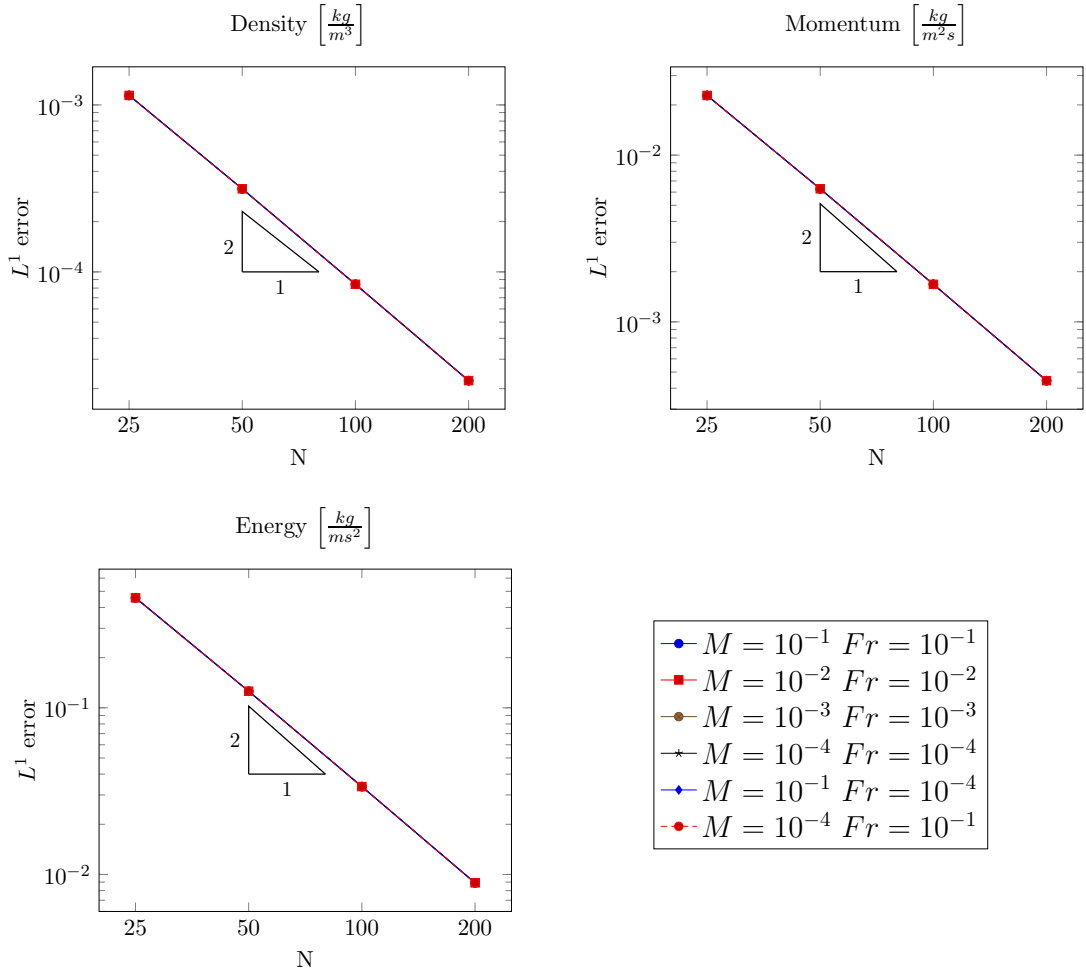

To transform the initial data \eqref{eq:ExactInitDim} into non-dimensional quantities, we define 
the following reference values
\begin{equation*}
x_r = 1m, ~u_r = 1\frac{m}{s}, ~\rho_r = 1 \frac{kg}{m^3},~ p_r = \frac{1}{M^2}\frac{kg}{m s^2}, ~ \Phi_r = \frac{1}{Fr^2} \frac{m^2}{s^2}.
\end{equation*} 
We use different values for $M$ and $Fr$ to show that our scheme is second order accurate independently of the chosen regime.
In the computations we use exact boundary conditions and $\gamma = 5/3$.
As can be seen from Table \ref{tab:Convergence2D} the error and the convergence rates are of the same magnitude for all displayed Mach numbers and we achieve the expected second order accuracy.
In addition, to illustrate that the accuracy is independent of the Mach number, we have plotted the $L^1$- error in Figure \ref{fig:error}. 
Due to the limiting procedure that we apply on the slopes in the reconstruction step to ensure the positivity property, we are not recovering a full second order convergence. 
Using unlimited slopes in the reconstruction step however will lead to the full second order.  
\begin{table}
	\begin{center}
		\small
		\renewcommand{\arraystretch}{1.25}
		\begin{tabular}{ccccccccccc}
			$M$&$Fr$&$N$ & $\rho \left[\frac{kg}{m^3}\right] $ & & $\rho u_1 \left[\frac{kg}{m^2 s}\right]$ & & $\rho u_2 \left[\frac{kg}{m^2s}\right]$ & & $E \left[\frac{kg}{m s^2}\right]$ & \\
			\hline \hline
			\multirow{4}{*}{$10^{-1}$}& \multirow{4}{*}{$10^{-1}$}& 25  &1.139E-003 &---&  2.278E-002  &---& 2.278E-002  &---& 4.562E-001&---\\
			&&50  &3.142E-004  &1.858& 6.276E-003 &1.859&  6.276E-003  &1.859& 1.257E-001&1.859\\
			&&100 & 8.427E-005  &1.898& 1.680E-003 &1.901& 1.680E-003 &1.901&  3.366E-002&1.901\\
			&&200 & 2.232E-005 &1.916&  4.438E-004  &1.920& 4.438E-004 &1.920&  8.894E-003&1.920\\
			\hline
			\multirow{4}{*}{$10^{-2}$}&\multirow{4}{*}{$10^{-2}$}&  25 &1.140E-003  &---& 2.280E-002 &---&  2.280E-002 &---&  4.567E-001  &--- \\
			&&50 & 3.144E-004  &1.859& 6.280E-003  &1.860& 6.280E-003 &1.860&  1.258E-001 &1.859 \\
			&&100 & 8.430E-005 &1.899&  1.680E-003 &1.901&  1.680E-003 &1.901&  3.367E-002  &1.901\\ 
			&&200 & 2.233E-005 &1.916&  4.441E-004  &1.919& 4.441E-004 &1.919&  8.901E-003&1.919\\
			\hline 
			\multirow{4}{*}{$10^{-3}$}&\multirow{4}{*}{$10^{-3}$} & 25 &1.141E-003  &---& 2.281E-002 &---&  2.281E-002 &---&  4.569E-001  &---\\
			&&50 & 3.144E-004  &1.859& 6.280E-003 &1.861&  6.280E-003 &1.861& 1.258E-001  &1.860\\ 
			&&100 & 8.431E-005 &1.898&  1.680E-003 &1.901&  1.680E-003 &1.901&  3.368E-002  &1.901\\ 
			&&200 & 2.233E-005  &1.916& 4.441E-004 &1.919&  4.441E-004  &1.919& 8.901E-003&1.919\\
			\hline
			\multirow{4}{*}{$10^{-4}$} &\multirow{4}{*}{$10^{-4}$}& 25 & 1.141E-003  &---& 2.280E-002 &---&  2.280E-002 &---&  4.582E-001 &---\\
			&&50 &3.143E-004 &1.860&  6.277E-003  &1.860& 6.277E-003 &1.860& 1.257E-001 &1.864\\ 
			&&100 &8.430E-005 &1.898&  1.680E-003 &1.901&  1.680E-003 &1.901&  3.367E-002 &1.901\\ 
			&&200 &2.233E-005 &1.916& 4.441E-004 &1.919& 4.441E-004 &1.919&  8.900E-003&1.919\\ 
			\hline 
			\multirow{4}{*}{$10^{-4}$} &\multirow{4}{*}{$10^{-1}$}& 25 & 1.141E-003 &---&  2.280E-002 &---&  2.280E-002 &---&  4.581E-001 &---\\
			&&50 & 3.143E-004  &1.860& 6.277E-003 &1.860& 6.277E-003 &1.860&  1.257E-001  &1.864\\ 
			&&100 & 8.430E-005 &1.898&  1.680E-003 &1.901&  1.680E-003 &1.901&  3.367E-002 &1.901 \\
			&&200 & 2.233E-005 &1.916&  4.441E-004 &1.919&  4.441E-004  &1.919& 8.900E-003 &1.919\\
			\hline
			\multirow{4}{*}{$10^{-1}$} &\multirow{4}{*}{$10^{-4}$}&25 &1.139E-003  &---& 2.278E-002 &---&  2.278E-002 &---&  4.562E-001 &---\\
			&&50 &3.142E-004 &1.858& 6.276E-003 &1.859&  6.276E-003  &1.859& 1.257E-001 &1.859\\ 
			&&100 &8.427E-005 &1.898&  1.680E-003 &1.901&  1.680E-003 &1.901&  3.366E-002 &1.901\\ 
			&&200 &2.232E-005 &1.916& 4.438E-004 & 1.920&  4.438E-004 & 1.920& 8.894E-003 &1.920\\
			\hline
		\end{tabular}
		\caption{$L^{1}$-error and convergence rates for different Mach and Froude numbers.}
		\label{tab:Convergence2D}
	\end{center}
\end{table}
\subsection{A stationary vortex in a gravitational field}

\begin{figure}[htpb]
	\centering
	\includegraphics[scale=0.65]{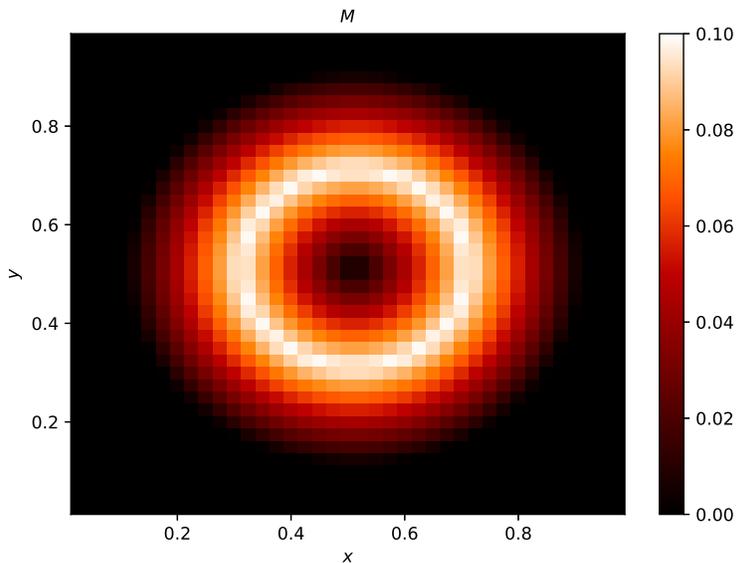}
	\caption{Initial Mach number distribution for $M=10^{-1}$.}
	\label{fig:GrafGreshoInit}
\end{figure}

\begin{figure}[htpb]
	\centering
	\includegraphics[scale=0.4]{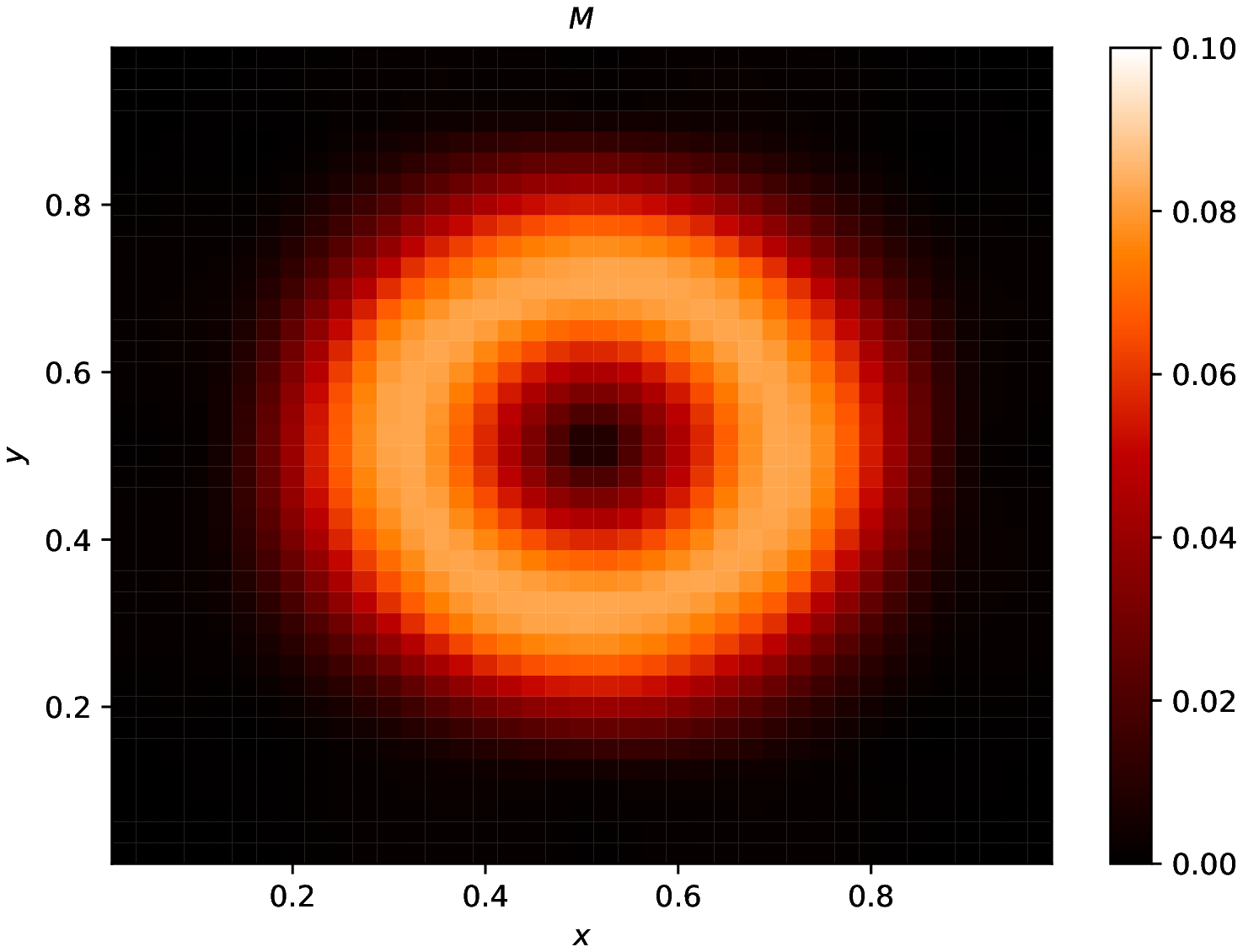}
	\includegraphics[scale=0.4]{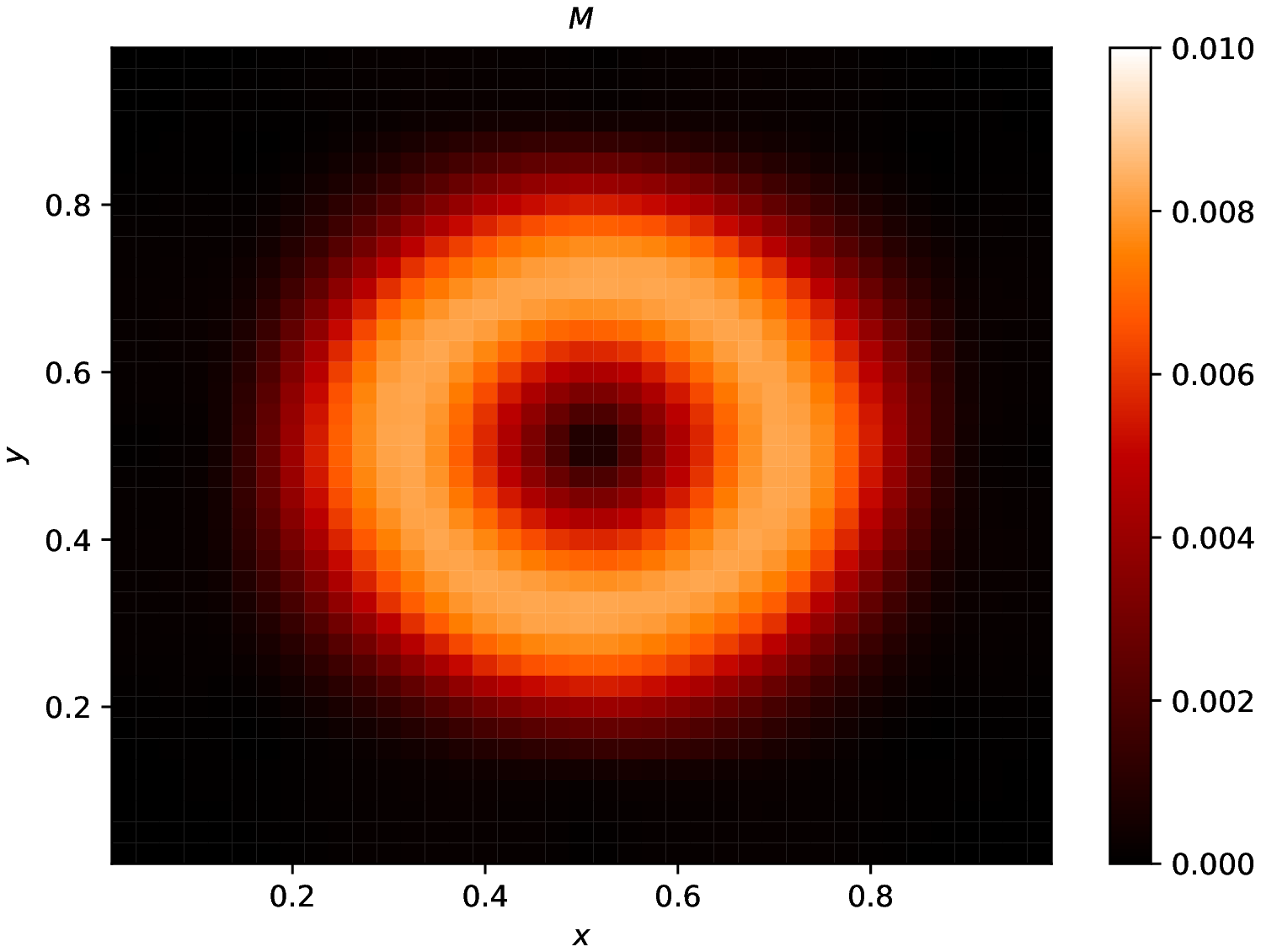}\\
	\includegraphics[scale=0.4]{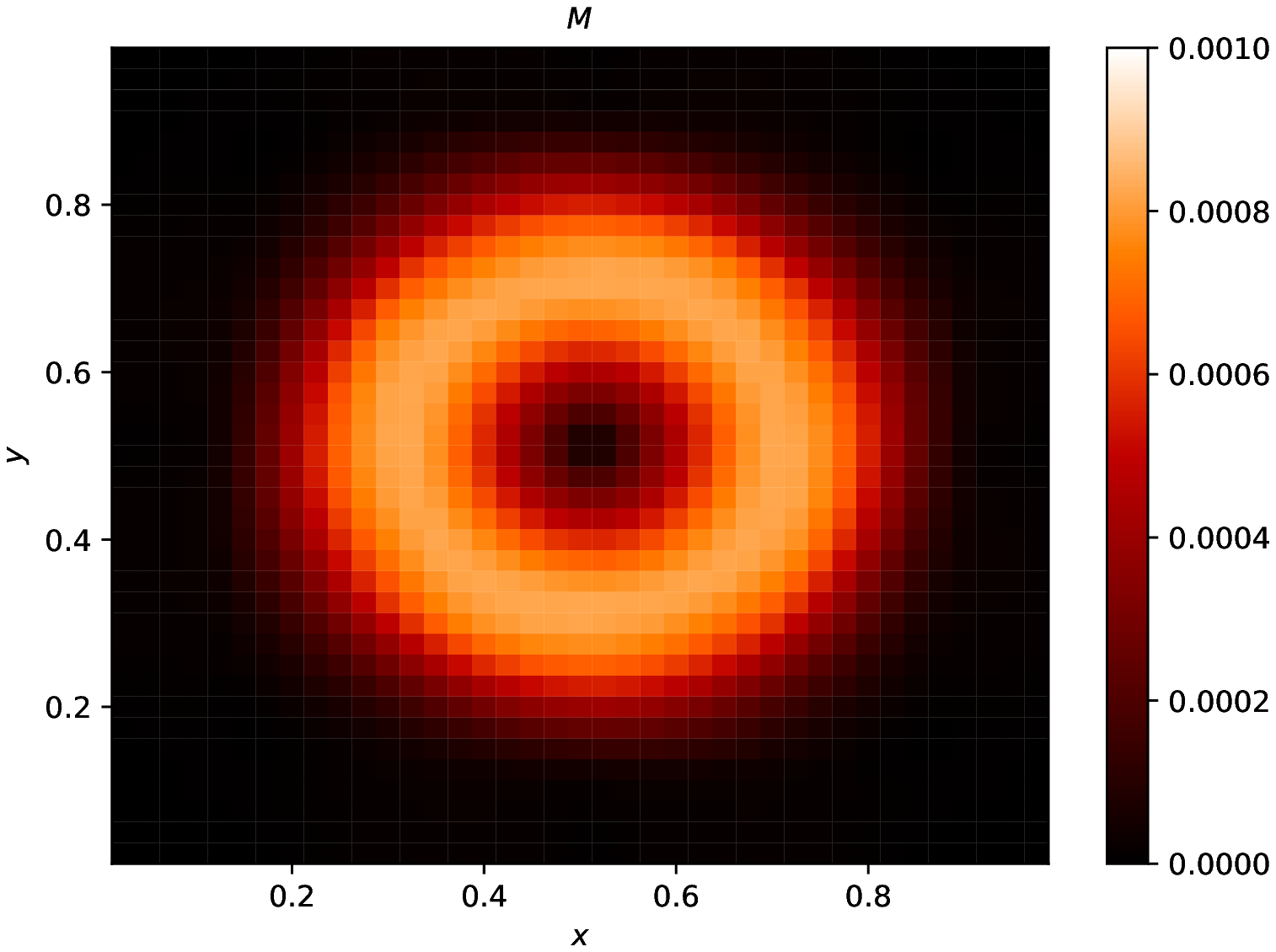}
	\includegraphics[scale=0.4]{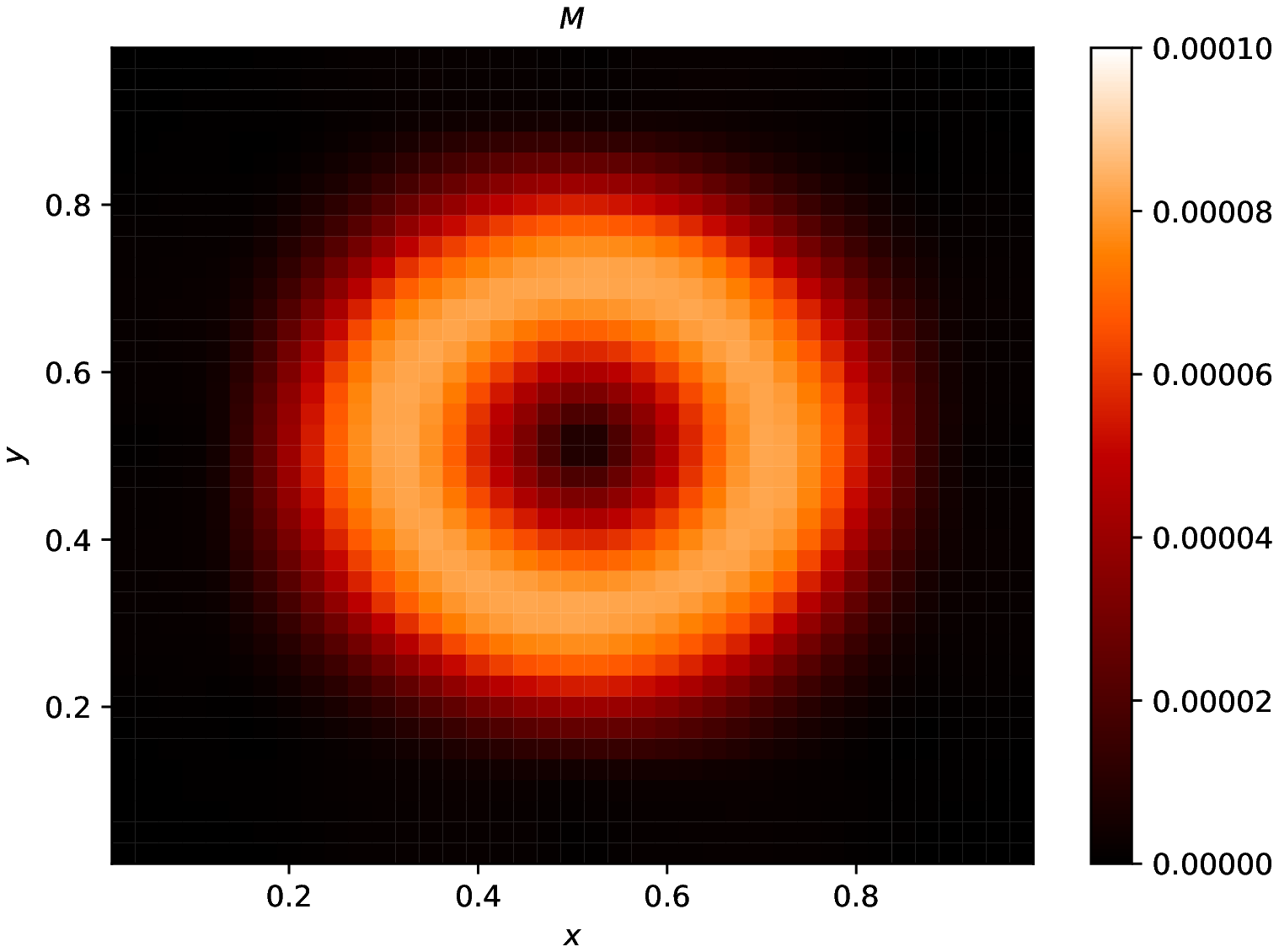}
	\caption{Mach number distribution for different maximal Mach numbers at $t=1$.\\ Top left: $M=10^{-1}$. Top right: $M=10^{-2}$, bottom left: $M=10^{-3}$, bottom right: $M=10^{-4}$}
	\label{fig:GrafGresho}
\end{figure}

\begin{figure}[htpb]
	\centering
	\includegraphics[scale=0.65]{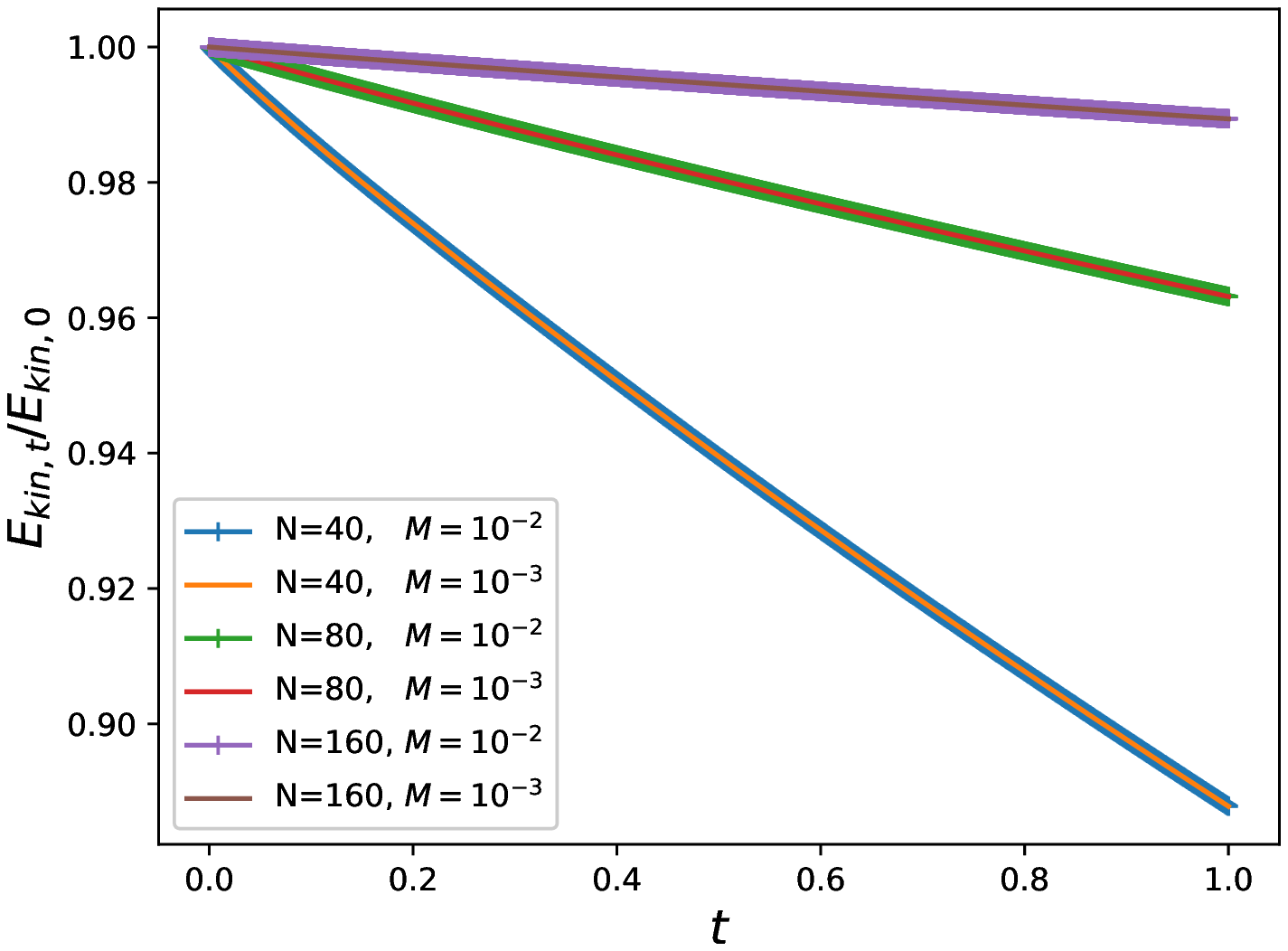}
	\caption{Loss of kinetic energy for different grids and Mach numbers after one full turn of the vortex (non-dimensional).}
	\label{fig:EkinGrafGresho}
\end{figure}

With this test-case, we want to demonstrate the low Mach properties of our scheme.
For the derivation of a vortex in a gravitational field, we follow the derivation of the Gresho vortex test case for the homogeneous Euler equations \cite{MiczekRoepkeEdelmann2015}.
It fulfils the divergence free property $\nabla \cdot \mathbf{u} = 0$ and the orthogonality property $\mathbf{u} \cdot \nabla \Phi = 0$ of the well-prepared data $\Omega_{wp}$.
To derive the vortex, we consider the non-dimensional Euler equations \eqref{sys:Euler_LMG} in radial coordinates $(r,\theta)$.
The vortex is constructed such that it is axisymmetric, stationary and has zero radial velocity.
A solution has to satisfy 
\begin{equation*}
\frac{1}{M^2}\partial_r p = \frac{\rho u_\theta^2}{r} - \rho \frac{\partial_r \Phi}{Fr^2},
\label{eq:StationaryEQGresho}
\end{equation*}
where $u_\theta$ is the angular velocity.
The pressure is split into a hydrostatic pressure $p_0$ and a pressure $p_2$ associated with the centrifugal forces and in total is given by $p = p_0 + M^2 p_2$ and has to satisfy
\begin{equation*}
\partial_r p_0 = -\frac{M^2}{Fr^2} \rho \partial_r \Phi, ~~~ \partial_r p_2 = \rho \frac{u_\theta^2(r)}{r}.
\end{equation*}
We choose an isothermal hydrostatic pressure $p_0 = R T \rho$ and the density is given according to \eqref{eq:isoth} by
\begin{equation*}
\rho = \exp\left(-\frac{M^2}{Fr^2}\frac{\Phi}{R T}\right).
\end{equation*}
The pressure $p_2$ is then given as
\begin{equation}
\label{eq:p2_int}
p_2 = \int_0^r \exp\left(-\frac{M^2}{Fr^2}\frac{\Phi(s)}{\chi}\right)\frac{u_\theta(s)^2}{s} ds.
\end{equation}
The velocity profile $u_\theta$ is defined piecewise as in the Gresho vortex test case as
\begin{equation*}
u_\theta(r) = \frac{1}{u_r}
\begin{cases}
5 r &\text{ if } ~~~r \leq 0.2, \\
2 - 5 r &\text{ if } ~~~0.2 < r \leq 0.4, \\
0 &\text{ if } ~~~r > 0.4.
\end{cases}
\end{equation*}
To fully determine $p_2$ a continuously differentiable gravitational potential has to be given. 
We define it piecewise as 
\begin{align*}
\Phi(r) = 
\begin{cases}
12.5 r^2 &\text{ if } ~~~ r \leq 0.2 \\
0.5 - \ln(0.2) + \ln(r) &\text{ if } ~~~ 0.2 < r \leq 0.4 \\
\ln(2) - 0.5\frac{r_c}{r_c - 0.4} + 2.5\frac{r_c}{r_c - 0.4}r - 1.25\frac{1}{r_c - 0.4} r^2 &\text{ if } ~~~0.4 < r \leq r_c \\
\ln(2) - 0.5 \frac{r_c}{r_c - 0.4} + 1.25 \frac{r_c^2}{r_c - 0.4}&\text{ if } ~~~r > r_c
\end{cases}.
\end{align*}
This choice of $\Phi$ ensures the use of periodic boundary conditions since $\Phi$ is constant at the boundary and thus we can simulate a closed system.
Then we can compute the pressure $p_2$ according to \eqref{eq:p2_int} and it is piecewise defined as 
\begin{align*}
p_2(r) = \frac{Fr^2 R T}{M^2 ~u_r^2}
\begin{cases}
p_{21}(r) &\text{ if } ~~~r \leq 0.2 \\
p_{21}(0.2) + p_{22}(r) & \text{ if } ~~~ 0.2 < r \leq 0.4 \\
p_{21}(0.2) + p_{22}(0.4) &\text{ if } ~~~ r > 0.4
\end{cases}
\end{align*}
with 
\begin{align*}
\begin{split}
p_{21}(r) = &\left(1 - \exp\left(-12.5\frac{ M^2}{Fr^2 RT} r^2\right)\right),\\
p_{22}(r) = &\frac{1}{\left(Fr^2 RT-M^2\right) \left(Fr^2 RT-0.5 M^2\right)}\exp\left(\frac{(-0.5 + \ln(0.2)) M^2}{Fr^2 RT}\right)\\ 
&\left(r^{-\frac{M^2}{Fr^2 RT}}\right.
\left(M^4 (r (10-12.5 r)-2)-4 Fr^4 \chi^2+Fr^2 M^2 (r (12.5 r-20)+6) RT\right)\\
&\left.+\exp\left(\frac{-\ln(0.2) M^2}{Fr^2 RT}\right) \left(4 Fr^4 RT^2-2.5 Fr^2 M^2 RT+0.5 M^4\right)\right).
\end{split}
\end{align*}
The reference values are defined as $x_r = 1 m$, $\rho_{r} = 1 \frac{kg}{m^3}$,  $u_{r} = 2 \cdot 0.2~\pi\frac{m}{s}$, $t_r=1\frac{m}{u_r}$ and $RT = \frac{1}{M^2} \frac{m^2}{s^2}$. 
The computations are carried out with $\gamma = 5/3$ and $M = Fr$ on the domain $D = [0,1]^2$.
In Figure \ref{fig:GrafGreshoInit} the initial Mach number distribution for the vortex for $M = 0.1$ is given.
In Figure \ref{fig:GrafGresho}, the Mach number distribution for different maximum Mach numbers are compared for $N = 40$ at $t=1$ which corresponds to one turn of the vortex. 
We see that the accuracy of the vortices are comparable independently of the chosen Mach number and they show the same amount of diffusiveness despite of the coarse grid used.
The periodic boundary conditions allow us to model a closed system and we can monitor the loss of kinetic energy during the simulation which is depicted in Figure \ref{fig:EkinGrafGresho}.
The graphs for the Mach numbers $M=10^{-2}$ and $M=10^{-3}$ are superposed which shows that the loss of kinetic energy is independent of the Mach number. 
This is in agreement with the theoretical results and demonstrates the low Mach number properties of the scheme. 
We remark that although using the second order scheme, we do not expect to get second order convergence due to the lack of smoothness in the velocity profile $u_\theta$ and therefore also in the energy.
\subsection{Rising bubble test case}
This test case is taken from \cite{MendezNCarroll1994} and models a rising bubble which has a higher temperature than the background atmosphere on the domain $D = [0 km, 10 km]\times [0 km, 15 km]$.
The gravitation acts along the $y$-direction and is given by 
\begin{equation*}
\Phi(x,y) = g y ~\frac{m^2}{s^2},
\end{equation*}
where $g= 9.81\frac{m}{s^2}$ is the gravitational acceleration.
The stratification of the atmosphere is given in terms of the potential temperature $\theta$ defined by 
\begin{equation*}
\theta = T \left(\frac{p_0}{p}\right)^{\frac{R}{c_p}},
\end{equation*}
where $c_p$ is the specific heat at constant pressure and $p_0 = 10^5 ~\frac{kg}{ms^2}$, denotes a reference pressure taken at sea level.
Pressure, potential temperature and density are connected by the following relation
\begin{equation}
\label{eq:RisingBubbleEOS}
p = p_0 \left(\frac{\theta R}{p_0}\right)^\gamma \rho^\gamma = \chi \rho^\gamma,
\end{equation} 
where $c_v$ is the specific heat at constant volume and $R = c_p - c_v$.
Comparing \eqref{eq:RisingBubbleEOS} to \eqref{eq:EOS}, the atmosphere is isentropic with the polytropic coefficient $\Gamma = \gamma$.
We set $p(x,0) = p_0$ and $\theta = 300K$. 
Therefore we have \begin{equation*}
\rho(x,0) = \frac{p_0}{\theta R}
\end{equation*}
and the hydrostatic equilibrium is given by \eqref{eq:poly}. To transform the data into non-dimensional quantities, we define the following reference values
\begin{equation*}
x_r = 10000~m, ~ t_r = 10000~s, ~u_r = 1\frac{m}{s}, ~\rho_r = 1 \frac{kg}{m^3}.
\end{equation*}
The scaling of the remaining variables is given in Table \ref{tab:Units+Scaling}.

The bubble is modelled as a disturbance in the potential temperature centred at $(x_c,y_c) = (5km,2.75km)$ as 
\begin{equation*}
\Delta \theta = 
\begin{cases}
\Delta \theta_0 \cos^2\left(\frac{\pi r}{2}\right) & \text{ if } r \leq 1\\
0 & \text{ else}
\end{cases}
\end{equation*}
where $\Delta \theta_0 = 6.6K$
and 
\begin{equation*}
r = \left(\frac{x - x_c}{r_0}\right)^2 + \left(\frac{y - y_c}{r_0}\right)^2
\end{equation*}
with the factor $r_0 = 2.0 km$.
The resulting perturbation in the pressure can be calculated from equation \eqref{eq:RisingBubbleEOS}.

\begin{figure}[htpb]
	\centering
	\includegraphics[scale=0.45]{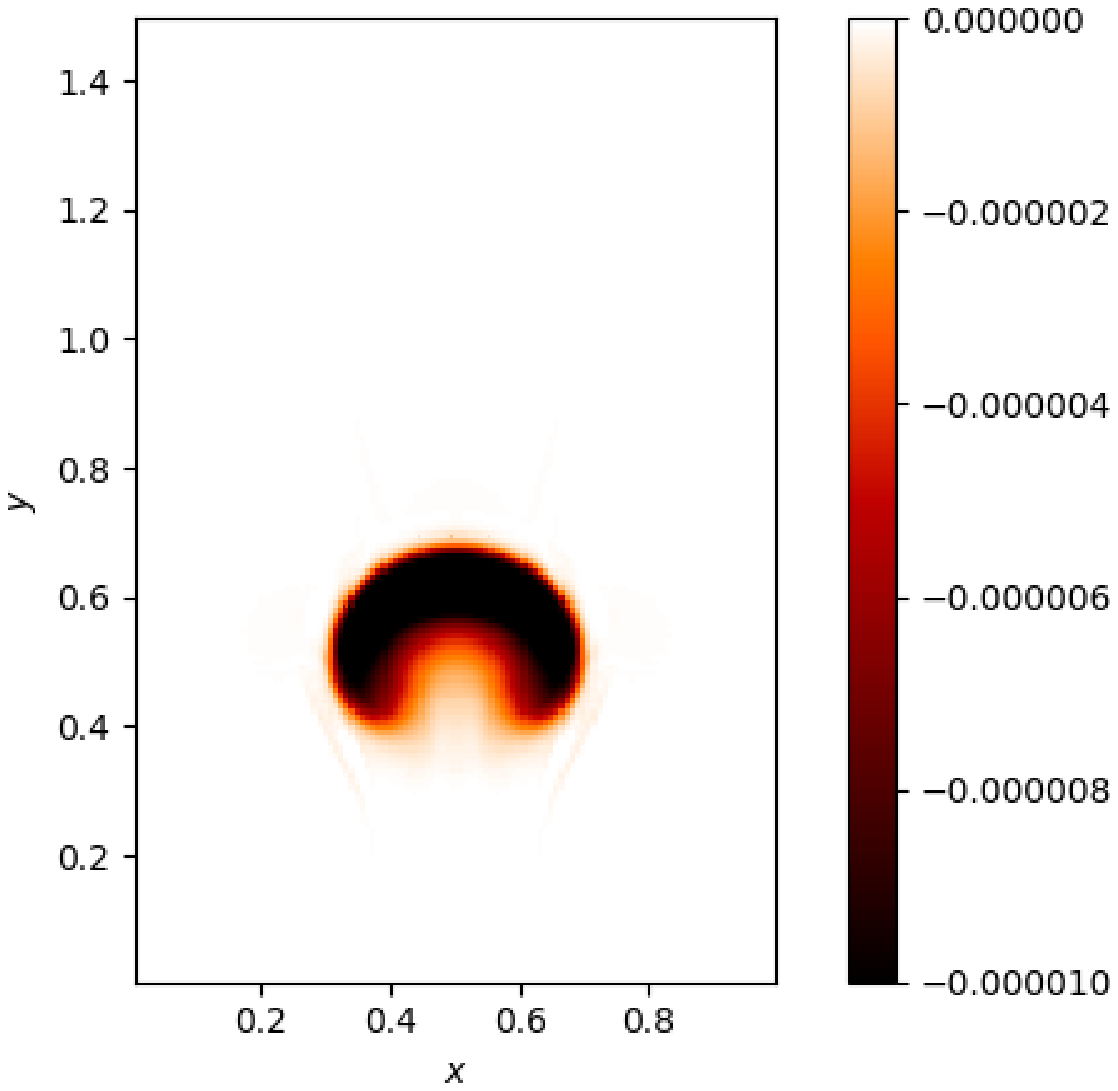}
	\includegraphics[scale=0.45]{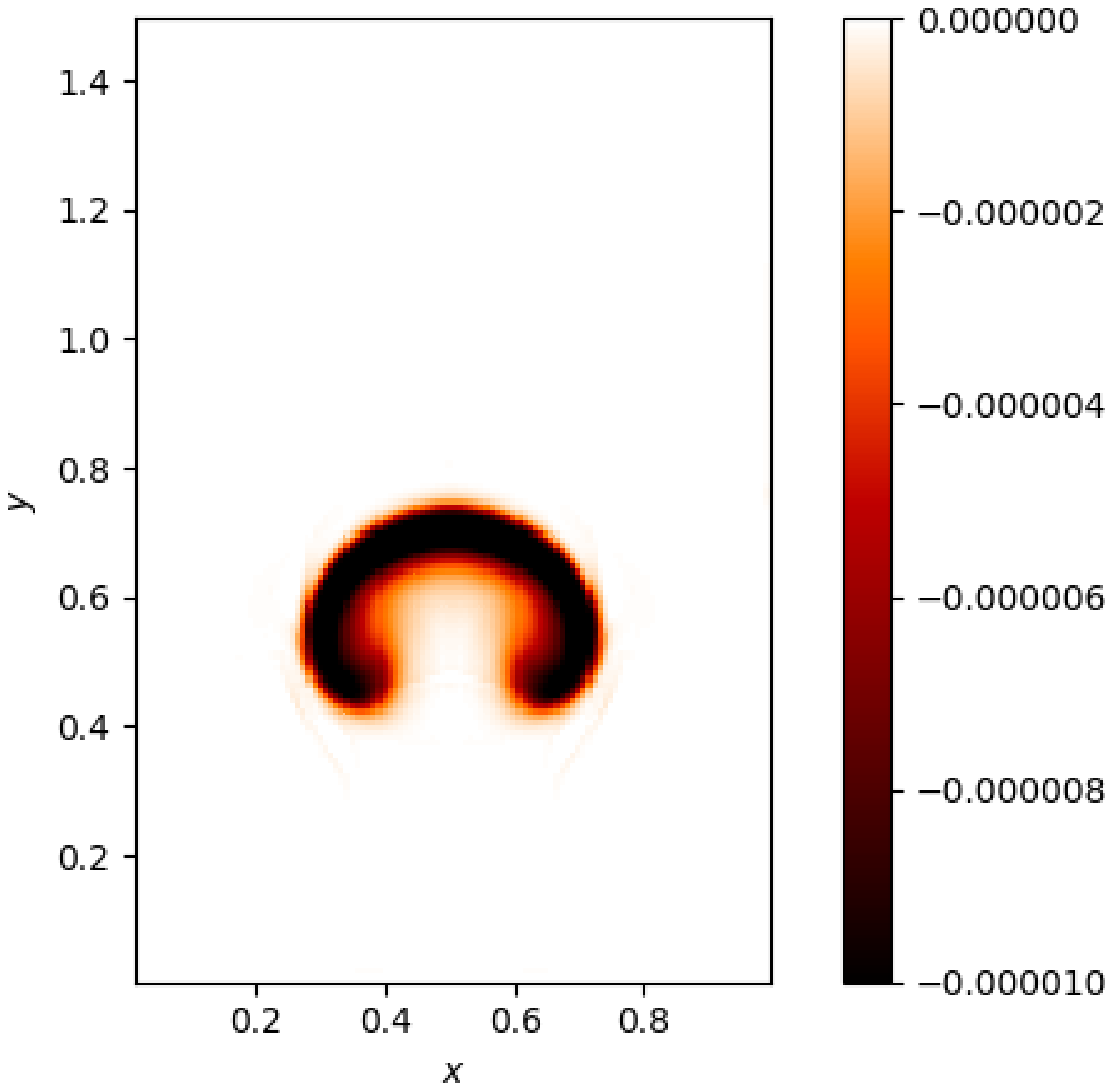}\\
	\includegraphics[scale=0.45]{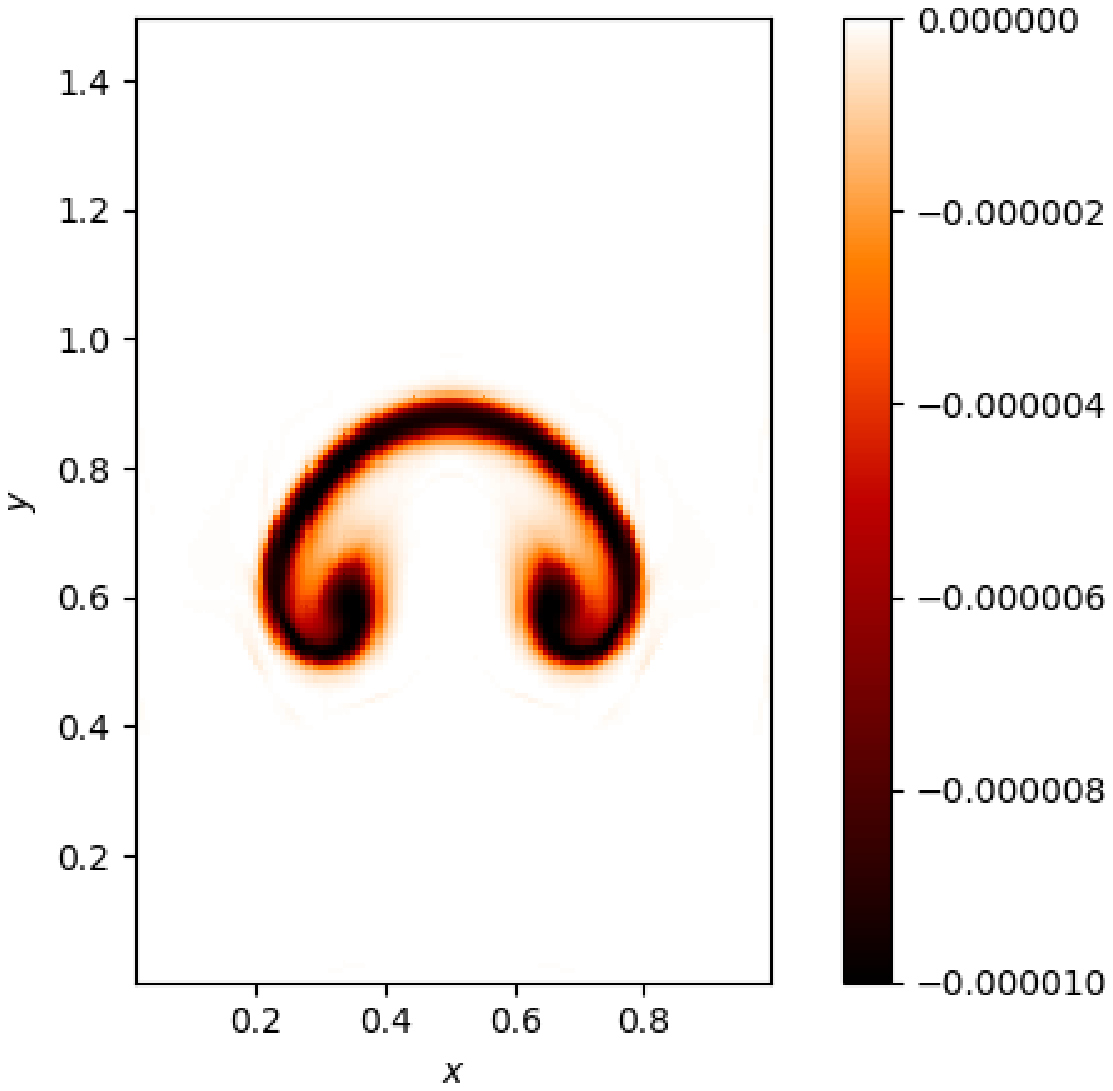}
	\includegraphics[scale=0.45]{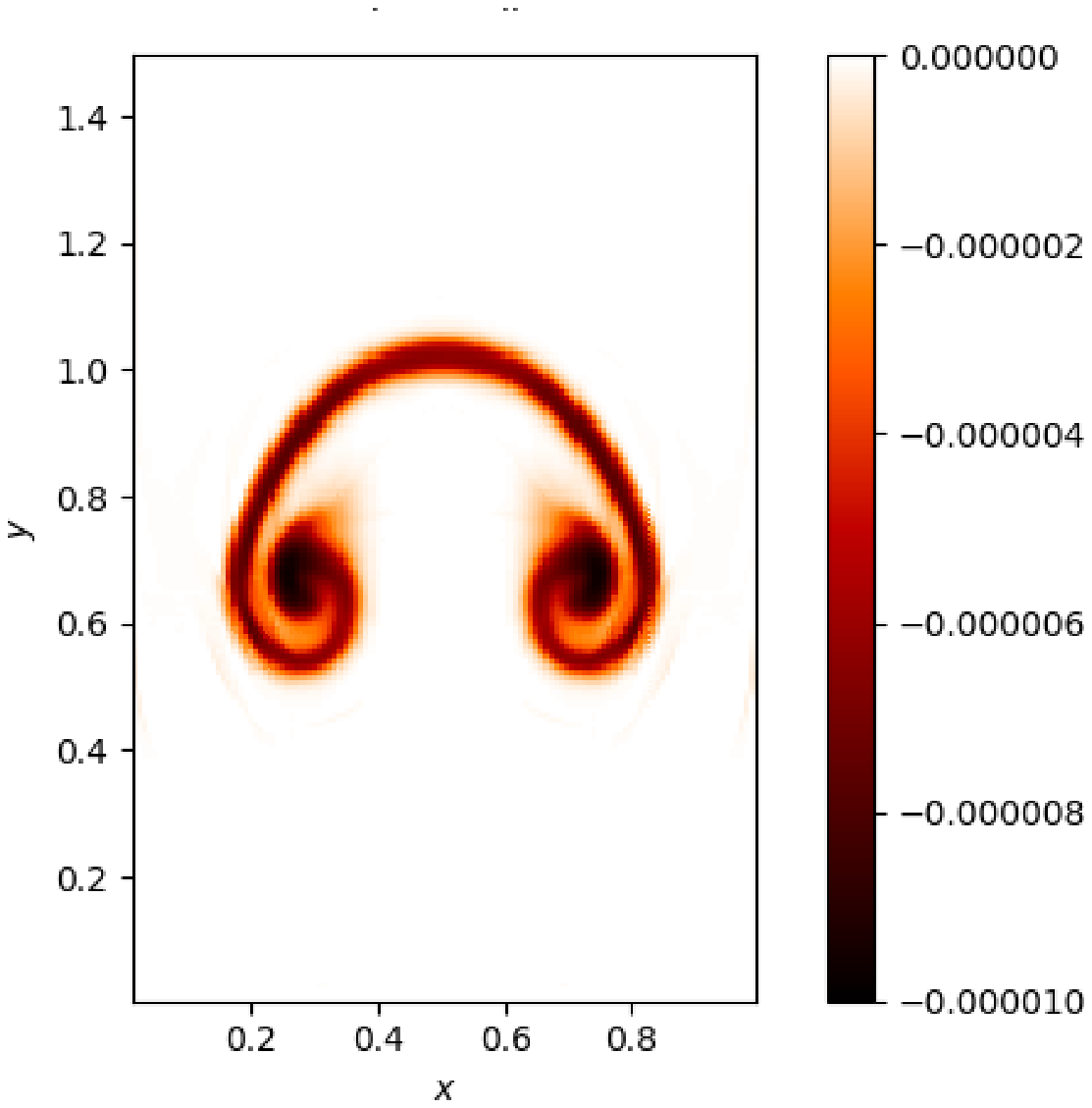}
	\caption{Density perturbation from the rising bubble test case from top right to bottom left at times $t=0.07, 0.09, 0.13, 0.18$.}
	\label{fig:RisingBubble_wb_bc}
\end{figure}

In the simulation, we choose $\gamma = 1.4$ as it is modelled air as a diatomic gas with the corresponding specific gas constant $R_s = 287.058 \frac{m^2}{s^2 K}$.
This setting results in a reference Mach number of $M=10^{-2}$ and we chose $Fr = M$.
In Figure \ref{fig:RisingBubble_wb_bc}, we show the density perturbation at different times $t$. 
It is computed with the second order scheme on a grid of $120$ cells in x-direction and $180$ cells in y-direction which results into a uniform space discretization. 
At the boundaries, we have imposed the background atmosphere. 
\begin{table}
	\centering
	\renewcommand{\arraystretch}{1.5}
	\begin{tabular}{c c c}
		\text{quantity } & \text{SI unit} & \text{ scaling }\\
		\hline \hline
		$ x $ & $\left[m\right]$ & $x_r$ \\ 
		$ t $ & $\left[s\right]$ & $t_r$ \\
		$\rho$ & $\left[\frac{kg}{m^3}\right]$ & $\rho_r$ \\
		$u, c$ & $\left[\frac{m}{s}\right]$ & $u_r = \frac{x_r}{t_r}$, $M = \frac{u_r}{c_r}$\\
		%$ g $ & $\left[\frac{m}{s^2}\right]$ & $g_r$ \\
		$ p  $ & $\left[\frac{kg}{m ~s^2}\right]$ & $p_r = R_s \rho_r \theta_r$, $p_r = \rho_r c_r^2$\\
		$ \Phi$ & $\left[\frac{m^2}{s^2}\right]$ & $\Phi_r = \frac{u^2_r}{Fr^2}$ \\%= x_r g_r$ \\
		$ R_s$ & $\left[\frac{m^2}{s^2 K} \right]$ & ---\\
		$ T,\theta$ & $\left[K\right]$ & $\theta_r = \frac{u_r^2}{R_s ~M^2}$\\
		\hline
	\end{tabular}
	\caption{Overview over units and scaling relations of the physical quantities used in the test cases in Section 6.}
	\label{tab:Units+Scaling}
\end{table}

\section{Conclusion} 
We have extended the second order all-speed IMEX scheme given in \cite{ThomannZenkPuppoKB2019} developed for the homogeneous Euler equations to treat a gravitational source term.
It is done in such a way that the new scheme inherits the positivity preserving property of the density and internal energy, as well as the scale independent diffusion and the AP property. 
In addition it is well-balanced for arbitrary hydrostatic equilibria. 
To show the AP property of the new IMEX scheme, we have defined a set of well-prepared data that consists to leading order of the hydrostatic equilibria where the velocity is divergence free and orthogonal to the direction of the gravitational potential. 
The resulting limit equations are the incompressible Euler equations with a gravitational source term. 
To numerically verify the low Mach properties of our scheme, we have developed a stationary vortex in a gravitational field which is well-prepared. 
With the help of this new test case we can demonstrate the scale independent diffusion of our scheme as it is already standard for the homogeneous case. 
The numerical results are concluded with a rising bubble test case to illustrate the applicability of our scheme.

\section*{Aknowledgements}
G. Puppo ackwowledges the support by the GNCS-INDAM 2019 research project and A. Thomann the support of the INDAM-DP-COFUND-2015, grant number 713485.
The authors would like to thank Markus Zenk for fruitful discussions and useful comments and suggestions. 

\bibliographystyle{unsrt}
\bibliography{lit.bib}

\end{document}